\documentclass[12pt]{article}
\usepackage{color}
\usepackage{epsfig,psfrag,amsmath,amssymb,latexsym}
\usepackage{amscd}
\usepackage{amsfonts}
\usepackage{graphicx}
\newtheorem{theorem}{Theorem}[section]

\newtheorem{corollary}{Corollary}[section]

\newtheorem{definition}{Definition}[section]
\newtheorem{example}{Example}[section]

\newtheorem{lemma}{Lemma}[section]

\newtheorem{proposition}{Proposition}[section]
\newtheorem{remark}{Remark}[section]

\newcommand{\NN}{\mathbb{N}}
\newcommand{\RR}{\mathbb{R}}

\newcommand{\PP}{\mathbb{P}}

\newcommand{\EE}{\mathbb{E}}
\newcommand{\II}{\mathbb{I}}

\newcommand{\FF}{\mathbb{F}}
\newcommand{\GG}{\mathbb{G}}
\newcommand{\M}{{\cal M}^{-1}}
\newcommand{\MP}{\mathcal{P}}
\newcommand{\p}{{\cal R}}
\newcommand{\q}{{\cal Q}}
\newcommand{\<}{\preccurlyeq}
\newcommand{\1}{{\bf 1}}
\newcommand{\0}{{\bf 0}}
\numberwithin{equation}{section}

\begin{document}

\title{\textbf{Hadamard functions of inverse $M$-Matrices}}

\author{Claude DELLACHERIE\thanks{Laboratoire Rapha\"el Salem,
UMR 6085, Universit\'e de Rouen, Site Colbert, 76821 Mont Saint
Aignan Cedex, France; Email: Claude.Dellacherie@univ-rouen.fr.
This author thanks support from Nucleus Millennium P04-069-F for
his visit to CMM-DIM at Santiago} , Servet
MARTINEZ\thanks{CMM-DIM; Universidad de Chile; Casilla 170-3
Correo 3 Santiago; Chile; Email: smartine@dim.uchile.cl. The
author's research is supported by Nucleus Millennium Information
and Randomness P04-069-F.} , Jaime SAN MARTIN\thanks{CMM-DIM;
Universidad de Chile; Casilla 170-3 Correo 3 Santiago; Chile;
Email: jsanmart@dim.uchile.cl. The author's research is supported
by Nucleus Millennium Information and Randomness P04-069-F.} . }
%EndAName

\maketitle

\begin{abstract} We prove that the class of GUM matrices is the
largest class of bi-potential matrices stable under Hadamard increasing
functions. We also show that any power $\alpha\ge 1$, in the sense of Hadamard functions,
of an inverse $M$-matrix is also inverse $M$-matrix showing a conjecture stated
in Neumann \cite{neummann1998}. We study the class of
filtered matrices, which include naturally the GUM matrices, and
present some sufficient conditions for a filtered matrix  to
be a bi-potential.
\end{abstract}



\smallskip
\noindent{\bf Running title:} Hadamard functions of inverse $M$-matrices.

\medskip

\section {Introduction and Basic Notations}

A nonnegative matrix $U$ is said to be a {\bf potential} if it is
nonsingular and its inverse satisfies
$$
\begin{array}{ll} &U^{-1}_{ij}\le 0 \hbox{ for } i\neq j ,\; U^{-1}_{ii} > 0
\cr  &\cr &\forall i\; \sum\limits_j U^{-1}_{ij} \ge 0,
\end{array}
$$
that is $U^{-1}$ is an $M$-matrix which is row diagonally dominant.
We denote this class of matrices by $\MP$. In addition we say that $U$ is a {\bf bi-potential}
if $U^{-1}$ is also column diagonally dominant. We denote this class by $bi\MP$.
We note that $\MP, bi\MP$ are contained in $\M$ the class of inverses of
$M$-matrices.
\medskip

The class of potentials matrices play an important
role in probability theory. They
represent the potential (from where we have taken the name) of a
transient continuous time Markov Chain $(X_t)_{t\ge 0}$, with
generator $-U^{-1}$. That is,
$$
U_{ij}=\int_0^\infty (e^{-U^{-1}t})_{ij}\; dt=\int_0^\infty
\PP_i\{X_t=j\} dt
$$
is the mean expected time expended at site $j$ when starting the
chain at site $i$. Clearly $U$ is a bi-potential if both $U$ and $U'$ are potentials.

\smallskip

To get a discrete time interpretation take
$K_0=\max\limits_i\{U^{-1}_{ii}\}$. For any $k\ge K_0$ the matrix
$P_k=\II-\frac{1}{k}U^{-1}$ is nonnegative, sub-stochastic
and
$$
U^{-1}=k(\II-P_k),
$$
If we can take $k=1$, then $ U^{-1}=\II-P$ (with $P=P_1$) and $U$
is the mean expected number of visits of a Markov chain
$(Y_n)_{n\in \NN}$ whose transition probability is given by $P$,
in fact
$$
U_{ij}=\sum\limits_{n\ge 0} P^n_{ij}=\sum_{n\ge 0} \PP_i\{Y_n=j\}.
$$

We notice that if $U$ is a potential then it is diagonally dominant on each row, in the sense
that for all $i,j$ we have $U_{ii}\ge U_{ji}$. The probabilistic prove of this fact
is that
$$
U_{ji}=f_{ji} U_{ii},
$$
where $f_{ji}\le 1$ is the probability that the Markov process $(X_t)$, starting from $j$ ever reaches
the state $i$. If $U$ is a bi-potential, then it is also diagonally dominant on each column. In this case
we just say it is diagonally dominant.

\bigskip

For any nonnegative matrix $U$ we define the following quantity
$$
\tau (U)=\inf\{t\ge 0:\; \II+tU \notin bi\MP\},
$$
which is invariant under permutations  that is $\tau(U)=\tau(\Pi
U\Pi')$. We point out that if $U$ is a positive matrix then
$\tau(U)>0$. We shall study some properties of this function
$\tau$. In particular we are interested in matrices for which
$\tau(U)=\infty$. The next result shows that $\tau(U)=\infty$ is a
generalization of the class $bi\MP$.

\begin{proposition}
\label{prop1} Assume $U$ is a nonnegative matrix, which is
nonsingular and $\tau(U)=\infty$, then $U \in bi\MP$.
\end{proposition}

\begin{proof} It is direct from the observation that
$$
t(\II+tU)^{-1}\mathop{\rightarrow}\limits_{t\to \infty} U^{-1}.
$$
$\Box$
\end{proof}

\begin{remark}
We shall prove later on that the reciprocal is also true: if $U$
is in the class $bi\MP$, then $\tau(U)=\infty$.
\end{remark}

The following notion will play an important role in this article.

\begin{definition} Given a matrix $B$ we say that a vector $\mu$
is a right equilibrium potential if
$$
B\mu=\1,
$$
where $\1$ is the constant vector of ones. Similarly it is defined
the notion of a left equilibrium potential, which is the right
equilibrium potential for $B'$. When $B$ is nonsingular we denote
the unique right and left equilibrium potentials by $\mu_B$ and
$\nu_B$.

We denote by $\bar \mu=\1'\mu$ the total mass of $\mu$. In the
nonsingular case, it is not difficult to see that $\bar \nu=\bar
\mu$.
\end{definition}

\medskip

Notice that for a matrix $U\in bi\MP$ the right and left equilibrium
potentials are nonnegative. This is exactly the same as the fact
that the inverse is row and column diagonally dominant.

\bigskip

\begin{definition}
Constant Block Form ({\bf CBF}) matrices can be defined
recursively in the following way: given two CBF matrices $A,B$ of
sizes $p$ and $n-p$ respectively, and numbers $\alpha, \beta$ we
produce the new CBF matrix by
\begin{equation}
\label{CBF} U=\begin{pmatrix} A &\alpha \1_p\1_{n-p}'\cr \beta
\1_{n-p}\1_p' & B
\end{pmatrix},
\end{equation}
where the vector $\1_p$ is the vector of ones of size $p$. We also
say that $U$ is in increasing CBF if $\; min\{A,B\}\ge
min\{\alpha,\beta\}$.
\end{definition}

\bigskip

We recall the following two definitions introduced in
\cite{mcdonald1995} and \cite{nabben1995}, that generalize the
concept of ultrametric matrices introduced in \cite{martinez1994}
(see also \cite{nabben1994}).

\begin{definition} A nonnegative CBF matrix $U$ is in
Nested Block Form ({\bf NBF}) if in (\ref{CBF}) $A, B$ are NBF
matrices and
\begin{itemize}
\item{} $0\le \alpha\le \beta$;
\item{} $\min\{A_{ij},A_{ji}\}\ge \alpha$ and $\min\{B_{kl},B_{lk}\}\ge
\alpha$;
\item{} $\max\{A_{ij},A_{ji}\}\ge \beta$ and $\max\{B_{kl},B_{lk}\}\ge
\beta$.
\end{itemize}
\end{definition}

\begin{definition} A nonnegative matrix $U$ of size $n$, is said to be
Generalized Ultrametric Matrix (GUM) if it is diagonally dominant
that is for all $i,j$ it holds  $U_{ii}\ge \max\{U_{ij},U_{ji}\}$,
and $n\le 2$, or $n>2$ and every three distinct elements $i,j,k$
has a preferred element. Assume this element is $i$, which means
\begin{itemize}
\item $U_{ij}=U_{ik}$;
\item $U_{ji}=U_{ki}$;
\item $\min\{U_{jk},U_{kj}\}\ge \min\{U_{ji},U_{ij}\}$;
\item $\max\{U_{jk},U_{kj}\}\ge \max\{U_{ji},U_{ij}\}$.
\end{itemize}
\end{definition}
By definition the transpose of a GUM matrix is also GUM.
We note that an ultrametric matrix is a symmetric GUM. The study of the chain
associated to an ultrametric matrix was done in \cite{dell1996} and for a GUM
in \cite{dell2000}.

\bigskip

In the next Theorem we summarize the main results in
\cite{mcdonald1995} and \cite{nabben1995} concerning GUM matrices.

\begin{theorem} Let $U$ be a nonnegative matrix.
\begin{itemize}
\item{} $U$ is a GUM matrix iff it is  permutation similar to a NBF.
\item{} If $U$ is a GUM, then it is nonsingular iff it does not
contain a row of zeros and no two rows are the same.
\item{} If $U$ is a non singular GUM then $U \in bi\MP$.
\end{itemize}
\end{theorem}
It is clear that if $U$ is a GUM then $\II+tU$ is a nonsingular
GUM. In particular $\tau(U)=\infty$.

\medskip

We introduce a main object of this article.

\begin{definition} Given a function $f$ and a matrix $U$, the
matrix $f(U)$ is defined as $f(U)_{ij}=f(U_{ij})$. We shall say
that $f(U)$ is a Hadamard function of $U$.

Given two matrices $A,B$ of the same size we denote by $A\odot B$
the Hadamard product of them. So $(A\odot B)_{ij}=A_{ij}B_{ij}$.
\end{definition}

Given a vector $a$ we denote by $D_a$ the diagonal
matrix, whose diagonal is $a$. For example we have
$D_aD_b=D_a\odot D_b=D_{a\odot b}$.

\bigskip

The class of CBF matrices (and its permutations) is closed under
Hadamard functions. Similarly, the class of increasing CBF (and
its permutations) is closed under increasing Hadamard functions.

On the other hand the class of NBF, and therefore also the class
of GUM matrices, is stable under Hadamard nonnegative increasing
functions. We summarize this result in the following proposition.

\medskip

\begin{proposition} Assume $U$ is a GUM and $f:\RR_+ \to \RR_+$ is
an increasing function. Then $f(U)$ is a GUM. In particular
$\tau(f(U))=\infty$, and if $f(U)$ is nonsingular then $f(U) \in
bi\MP$. A sufficient condition for $f(U)$ to be nonsingular is that
$U$ is nonsingular and $f$ is strictly increasing.
\end{proposition}

\begin{proof} It is clear that $f(U)$ is a GUM matrix and
therefore $\tau(f(U))=\infty$. Then, from Proposition \ref{prop1}
we have that $f(U) \in bi\MP$ as long as it is nonsingular. If $U$
is nonsingular then it does not contain a row (or column) of zeros
and there are not two equal rows (or columns). This condition is
stable under strictly non-negative functions, so the result
follows.$\Box$
\end{proof}

\bigskip

One of our main results is a sort of reciprocal of the previous
one. We shall prove that if $\tau(f(U))=\infty$ for all increasing
nonnegative functions $f$, then $U$ must be a GUM (see Theorem
\ref{establecreciente->GUM}).

\bigskip

The last concept we need for our work is the following.

\medskip

\begin{definition}
\label{classtau} We say that a nonnegative matrix $U$ is in class
$\mathcal{T}$ if
$$
\tau(U)=\inf\{t>0:\; (\II+tU)^{-1}\1 \ngeqq 0 \hbox{ or }
\1'(\II+tU)^{-1} \ngeqq 0\},
$$
and $\; \II+\tau(U)\, U$ is nonsingular whenever $\tau(U)<\infty$.
\end{definition}

\bigskip

We shall prove that every nonnegative matrix $U$ which is a
permutation of an increasing CBF, is in class $\mathcal{T}$.

\bigskip

We remark here that our purpose is to study Hadamard functions of matrices and
not spectral functions of matrices, which are quite different concepts.
For spectral functions of matrices there are deep and beautiful results
for the same classes of matrices we consider here. See for example the work
of Bouleau \cite{bouleau1989} for filtered operators. For $M$ matrices
see the works of Varga \cite{varga1968}, Micchelli and Willoughby \cite{micchelli1979},
Ando \cite{ando1980},  Fiedler and Schneider \cite{fiedler1983}, and the recent work of
Bapat, Catral and Neummann \cite{neummann2005} for $M$-matrices and inverse $M$-matrices.

\bigskip

\section{Main Results}

\begin{theorem}
\label{potencialfU} Assume $U\in \MP$ and $f:\RR_+\to \RR_+$ is a
nonnegative strictly increasing convex function. Then $f(U)$ is
nonsingular and $det(f(U))>0$. Also $f(U)$ has a right
nonnegative equilibrium potential. Moreover if $f(0)=0$ we have
$M=U^{-1}f(U)$ is an $M$-matrix.
If $U\in bi\MP$ then $f(U)$ also has a left nonnegative equilibrium potential.
\end{theorem}

Note that $H=f(U)^{-1}$ is not necessarily a $Z$-matrix, that is
for some $i\neq j$ it can happen that $H_{ij}>0$, as the following
example will show. Therefore the existence of a nonnegative right
equilibrium potential, which is
$$
\forall i\;\; H_{ii}+\sum\limits_{j\neq i} H_{ij}\ge 0.
$$
does not imply that the inverse is row diagonally dominant, that
is
$$
\forall i\;\; H_{ii}\ge \sum\limits_{j\neq i} |H_{ij}|.
$$

\begin{example} Consider the matrix
$$
P=\begin{pmatrix} 0 &\frac{1}{2} &0\cr \frac{1}{2} &0
&\frac{1}{2}\cr 0 &\frac{1}{2} &0
\end{pmatrix}
$$
Then $U=(\II-P)^{-1} \in bi\MP$. Consider the strictly convex
function $f(x)=x^2-\cos(x)+1$. A numerical computation gives
$$
(f(U))^{-1}\approx\begin{pmatrix} \;\;\;0.3590 &-0.0975
&\;\;\;0.0027\cr -0.0975 &\;\;0.2372 &-0.0975\cr \;\;\;0.0027
&-0.0975 &\;\;\;0.3590
\end{pmatrix},
$$
which is not a $Z$-matrix.
\end{example}

\medskip

We denote by $U^{(\alpha)}$ the Hadamard transformation of $U$
under $f(x)=x^\alpha$. In particular $U^{(2)}=U\odot U$. One of
our main results is the stability of $\M$ under powers. This
solves a conjecture stated on \cite{neummann1998}.

\begin{theorem}
\label{power} Assume $U\in \M$ and  $~\alpha\ge 1$. Then
$U^{(\alpha)}\in \M$. If $U^{-1}\in \MP$ then $(U^{(\alpha)})^{-1}\in \MP$.
If $~U\in bi\MP$ then $U^{(\alpha)} \in bi\MP$.
\end{theorem}

\medskip

The previous result has the following probabilistic
interpretation. If $U$ is the potential of a transient continuous
time Markov process then $U^{(\alpha)}$ is also the potential of a
transient continuous time Markov process. In the next result we
show the same is true for a potential of a Markov chain. An interesting
open question is what is the relation between the Markov
chain associated to $U$ and the one associated to $U^{(\alpha)}$.

\begin{theorem}
\label{markovchain} Assume that $U^{-1}=\II-P$ where $P$ is a
submarkov kernel, that is $P\ge 0$, $P\1\le \1$. Then for all
$\alpha \ge 1$ there is a submarkov kernel $Q(\alpha)$ such that
$(U^{(\alpha)})^{-1}=\II-Q(\alpha)$. Moreover if $P'\1\le \1$ then
$Q(\alpha)'\1\le \1$.
\end{theorem}

\medskip

The next result establishes that the class of GUM matrices is the largest class of potentials
stable under increasing Hadamard functions.

\begin{theorem}
\label{establecreciente->GUM} Let $U$ be a nonnegative matrix such
that $\tau(f(U))=\infty$, for all increasing nonnegative functions
$f$. Then, $U$ must be a GUM.
\end{theorem}

\medskip

\medskip

\begin{example}
\label{ejemplo1} Given $a,b,c,d \in \RR_+$ consider the
non-singular matrix
$$
U=\begin{pmatrix} 1 &0 &0 &0\cr 0 &1 &0 &0\cr a &b &1 &0\cr c &d
&0 &1 \end{pmatrix}.
$$
For all increasing non-negative functions $f$ and all $t>0$:
$(\II+tf(U))^{-1}$ is an $M$-matrix, while $U$ is not a GUM.
Moreover, $U$ is not a permutation of an increasing CBF. This
shows that the last Theorem does not hold if we replace the class
$bi\MP$ by the class $\M$.
\end{example}

\begin{theorem}
\label{CBFconvex} Let $U\in bi\MP$ and $f:\RR_+\to \RR_+$ be a
strictly increasing convex function. $f(U)$ is in $bi\MP$ if and
only if $f(U)$ belongs to class $\mathcal{T}$.
\end{theorem}

\begin{theorem}
\label{CBFisTau} If $U$ is a nonnegative increasing CBF matrix
then $U$ is in class $\mathcal{T}$.
\end{theorem}

\bigskip
As a corollary of the two previous theorems we obtain the following important result.
\begin{theorem}
Assume that $U\in bi\MP$ is an increasing CBF matrix, and
$f:\RR_+\to \RR_+$ is a nonnegative strictly increasing convex
function. Then $f(U)\in bi\MP$.
\end{theorem}

\bigskip

\section{Proof of Theorems \ref{potencialfU}, \ref{power}, \ref{markovchain} and \ref{CBFconvex}  }

Let us start with a useful lemma.

\begin{lemma}
\label{M->Mt} Assume $U \in \M$. Then
for all $t\ge 0$ we have  $(\II+t U)\in \M$. Also if $U \in \MP$ so is
$(\II+t U)$ and its right equilibrium potential is strictly positive.
In particular if $U\in bi\MP$ then so is
$\II+tU$ and its equilibrium potentials are strictly positive.
Similarly, let  $0\le s<t$ and assume $\II+tU \in bi\MP$, then
$\II+sU\in bi\MP$ and its equilibrium potentials are strictly
positive.
\end{lemma}

\begin{proof} For some  $k>0$ large enough, $U^{-1}=k(\II-N)$
where $N\ge 0$ (and  $N\1\le \1$ for the row diagonally dominant
case). In what follows we can assume that $k=1$ (it is enough to
consider the matrix $kU$ instead of $U$).

From the equality $(\II-N)(\II+N+N^2+\cdots N^p)=\II-N^{p+1}$ we
get that
$$
\II+N+N^2+\cdots N^p= U (\II-N^{p+1})\le U,
$$
and we deduce that the series $\sum\limits_{l=1}^\infty N^l$ is
convergent and its limit is $ U$.

Consider now the matrix
$$
N_t=t\left(\left(\II-\frac{1}{1+t}N\right)^{-1}-\II\right)=t\sum\limits_{l=1}^\infty
\left(\frac{1}{1+t}\right)^l N^l.
$$
We have that $N_t\ge 0$ (and  $N_t\1\le \1$ whenever $N\1\le \1$).
Therefore the matrix $\II-N_t$ is an $M$-matrix (which is row
diagonally dominant when $M$ is so). On the other hand we have
$$
\II+tU=\II+t(\II-N)^{-1}=(t\II+\II-N)(\II-N)^{-1}=(1+t)\left(\II-\frac{1}{1+t}N\right)(\II-N)^{-1},
$$
from where we deduce that $\II+tU$ is nonsingular and its inverse
is
$$
\begin{array}{ll}
(\II+tU)^{-1}&=\frac{1}{1+t}(\II-N)\left(\II-\frac{1}{1+t}N\right)^{-1}\cr
&\cr &=
\frac{1}{1+t}\left(\left(\II-\frac{1}{1+t}N\right)^{-1}-N\left(\II-\frac{1}{1+t}N\right)^{-1}\right)\cr
& \cr &=\frac{1}{1+t}\left(\sum\limits_{l=0}^\infty (1+t)^{-l} N^l
-\sum\limits_{l=0}^\infty (1+t)^{-l} N^{l+1}\right)\cr & \cr &=
\frac{1}{1+t}(\II-N_t).
\end{array}
$$
The only thing left to prove is that $N_t \1 < \1$ in the row
diagonally dominant case, that is when $N\1\le \1$. For that it is
enough to prove that $N^l\1 < \1$ for large $l$. From the equality
$U=\sum\limits_{l=1}^\infty N^l$, we deduce that
$\sum\limits_{l=1}^\infty N^l \1<\infty$, and therefore $N^l\1$
tend to zero as $l\to \infty$. This proves the claim.

When $k$ is not $1$ we have the following equality
$(\II+tU)^{-1}=\frac{k}{t+k}(\II-\frac{t}{k}
\sum\limits_{l=1}^\infty (\frac{k}{t+k})^l N^l)$, where
$N=\II-\frac{1}{k}U^{-1}$.

\medskip

Finally, assume that $\II+tU \in bi\MP$. Hence $\II
+\beta(\II+t U) \in bi\MP$ for all $\beta\ge 0$. This implies that
$$
\II+\frac{\beta}{1+\beta} t\, U \in bi\MP.
$$
Now, it is enough to take $\beta\ge 0$ such that
$s=\frac{\beta}{1+\beta} t$.
$\Box$
\end{proof}

\bigskip

This Lemma has two immediate important consequences.

\begin{corollary}  If $U\in bi\MP$ then
$\tau(U)=\infty$.
\end{corollary}

\begin{corollary}  Let $U$ a nonnegative matrix, then
$$
\tau(U)=\sup\{t\ge 0:\; \II+tU \in bi\MP\}
$$
\end{corollary}
\begin{proof} It is clear that $\tau(U)\le \sup\{t\ge 0:\; \II+tU \in
bi\MP\}$. On the other hand if $\II+tU \in bi\MP$ then we get $\II +s
U \in bi\MP$ for all $0\le s \le t$. This fact and the definition of
$\tau(U)$, implies the result. $\Box$
\end{proof}

\bigskip

\begin{proof} {(\bf Theorem \ref{potencialfU})} We first assume that $f(0)=0$.
We have that $U^{-1}=k(\II-P)$, for some $k>0$ and  $P$  a
sub-stochastic matrix. Without loss of generality we can
assume $k=1$, because it is enough to consider $k U$ instead of
$U$ and $\tilde f(x)=f(x/k)$ instead of $f$.

Consider $M=(U^{-1}f(U))$.
Take $i\neq j$ and compute
$$
M_{ij}=(U^{-1}f(U))_{ij}=(1-p_{ii})f(U_{ij})-\sum\limits_{k\neq i}
p_{ik}f(U_{kj}).
$$
Since $1-p_{ii}-\sum\limits_{k\neq i} p_{ik}\ge 0$, which is
equivalent to $\sum\limits_{k} p_{ik}\le 1$, and $f$ is convex we
obtain
$$
\left(1-\sum\limits_{k} p_{ik}\right)f(0)+\sum\limits_{k}
p_{ik}f(U_{kj})\ge f\left(\sum\limits_{k}
p_{ik}U_{kj}\right)=f(U_{ij}).
$$
The last equality follows from the fact that $U^{-1}=\II-P$. This
shows that $M_{ij}\le 0$. Consider now a positive vector $x$ such
that $y'=x'U^{-1}>0$ (see \cite{HJ}, Theorem 2.5.3). Then
$$
x'M=x'U^{-1}f(U)=y'f(U)>0,
$$
which implies, by the same cited
theorem, that $M$ is a $M$-matrix. In particular $M$ is nonsingular and
$det(M)>0$. So, $f(U)$ is nonsingular and $det(f(U))>0$. Consider
now $\rho$ the right equilibrium potential of $f(U)$. We have
$$
M\rho=U^{-1}f(U)\rho=U^{-1}\1=\mu_U\ge 0,
$$
then $\rho=M^{-1}\mu_U\ge 0$, because $M^{-1}$ is a nonnegative
matrix. This means that $f(U)$ posseses a nonnegative right
equilibrium potential. Since $f(U)$ is non singular we also have a left
equilibrium potential, which we do not know if it is nonnegative.
Then the first part is proven under the
extra hypothesis that $f(0)=0$.

Assume now $a=f(0)>0$, and consider $g(x)=f(x)-a$, which is a
strictly increasing convex function. Obviously $f(U)=g(U)+a\1\1'$,
so
$$
\mu_{f(U)}=\frac{1}{1+a\bar \mu_{g(U)}}\mu_{g(U)}\ge 0,\;\;
\nu_{f(U)}=\frac{1}{1+a\bar \mu_{g(U)}}\nu_{g(U)},
$$
where $\bar\mu_{g(U))}=\1'\mu_{g(U))}>0$, and we have used the fact that
$\bar\mu_{g(U))}=\bar\nu_{g(U))}$. Thus $f(U)$ has
a nonnegative right equilibrium potential, and a left equilibrium potential.
We need to prove that $f(U)$
is nonsingular, and $det(f(U))>0$. This follows immediately from
the equality
$$
f(U)=g(U)(\II+a\mu_{g(U)}\1').
$$
Indeed we have
$$
f(U)^{-1}=g(U)^{-1}-\frac{a}{1+a\bar \mu_{g(U)}}\mu_{g(U)}
(\nu_{g(U)})' \hbox{ and } det(f(U))=det(g(U))(1+a\bar\mu_{g(U)}),
$$
from which the result is proven.

In the bi-potential case use  $U'$ instead of $U$ to obtain the
existence of a nonnegative left equilibrium potential for $f(U)$.

$\Box$
\end{proof}

\medskip

\begin{proof}({\bf Theorem \ref{CBFconvex}})
Using the same ideas as above we can assume that $f(0)=0$. Also we
have  $U^{-1}(\II+t f(U))=M_t$ is a $M$-matrix, for all $t\ge 0$.
Therefore $\II+t f(U)$ is nonsingular for all $t$, and we denote
by $\mu_t$ and $\nu_t$ the equilibrium potentials for $\II+t
f(U)$.

\medskip

Assume first that $f(U)$ is in class $\mathcal{T}$  which means
that
$$
\tau(f(U))=\min\{t> 0:\; \mu_t \ngeqq 0 \hbox{ or } \nu_t \ngeqq
0\}.
$$
We prove that for all $t\ge 0$, $\mu_t, \nu_t$ are nonnegative.
Since
$$
M_t \mu_t=U^{-1}\1=\mu_U,
$$
we obtain that $\mu_t=M_t^{-1}\mu_U\ge 0$, because $M_t^{-1}$ is a
nonnegative matrix. Thus, $\tau(f(U))=\infty$ and  $f(U)$ is
nonsingular. From Proposition \ref{prop1} we get $f(U)\in bi\MP$.

\smallskip

Reciprocally if $f(U)\in bi\MP$ then $\tau(f(U))=\infty$ and the
result follows. $\Box$
\end{proof}

\bigskip

\begin{lemma} \label{lm2} Assume that $U\in \MP$. Then any principal
square submatrix $A$ of $U$ is also in class $\MP$. The same is true if we replace
$\MP$ by $bi\MP$.
\end{lemma}

\begin{proof}
By induction and a suitable permutation is enough to prove the
result for $A$ the restriction of $U$ to $\{1,\dots,n-1\}\times
\{1,\dots,n-1\}$, where $n$ is the size of $U$. Assume that
$$
U=\begin{pmatrix} A &b \cr c' &d
\end{pmatrix} \hbox{ and }
U^{-1}=\begin{pmatrix} \Lambda &-\zeta \cr -\varrho' &\theta
\end{pmatrix}.
$$
Since $A^{-1}=\Lambda-\frac{1}{\theta}\zeta \varrho'$ we obtain
that the off diagonal elements of $A^{-1}$ are non-positive. It is
quite easy to see that the result will follow as soon as
$A^{-1}\1\ge 0$.

Since $U\in \MP$ we have that $\Lambda \1 - \zeta \ge 0$ and
$\theta \ge \varrho' \1$. Therefore,
$$
A^{-1}\1=\Lambda\1-\frac{1}{\theta}\zeta \varrho'\1=
\Lambda\1-\frac{\varrho'\1}{\theta}\;\zeta \ge \Lambda\1-\zeta \ge
0.
$$

$\Box$
\end{proof}

In what follows given a vector $a$ we denote by $D_a$ the diagonal
matrix, whose diagonal is $a$.

\begin{lemma}
\label{generalequilibrium} Assume $U\in bi\MP$ and $\alpha\ge 1$. If
$$
U=\begin{pmatrix} A &b \cr c' &d
\end{pmatrix},
$$
then there exists a nonnegative vector $\eta$ such that
$$
A^{(\alpha)} \eta=b^{(\alpha)}.
$$
\end{lemma}

\begin{proof} We first perturb the matrix $U$ to have a positive
matrix. Consider $\epsilon>0$ and the positive matrix
$U_\epsilon=U+\epsilon \1\1'$. It is direct to prove that
$$
U_\epsilon^{-1}=U^{-1}-\frac{\epsilon}{1+\epsilon \bar \mu_U}\mu_U
(\nu_U)',
$$
where $\bar \mu_U=\1'\mu_U$ is the total mass of $\mu_U$. Then
$U_\epsilon \in bi\MP$ and its equilibrium potentials  are given by
$$
\mu_{U_\epsilon}=\frac{1}{1+\epsilon \bar\mu_U} \mu_U,\;\;
\nu_{U_\epsilon}=\frac{1}{1+\epsilon \bar\nu_U} \nu_U.
$$
We decompose the inverse of $U_\epsilon$ as
$$
U_\epsilon^{-1}=\begin{pmatrix} \Lambda_\epsilon &\zeta_\epsilon
\cr \varrho_\epsilon' &\theta_\epsilon
\end{pmatrix},
$$
and we notice that $U_\epsilon \zeta_\epsilon+ \theta_\epsilon
b_\epsilon=0$ which implies that
$$
b_\epsilon=U_\epsilon \lambda_\epsilon,
$$
with $\lambda_\epsilon=-\frac{1}{\theta_\epsilon}\zeta_\epsilon\ge
0$. Also we mention here that $\lambda_\epsilon$ is a
sub-probability vector, that is $\1'\lambda_\epsilon\le 1$. This
follows from the fact that $U_\epsilon^{-1}$ is column diagonally
dominant.

Take now the matrix $V_\epsilon=D^{-1}_{b_\epsilon} U_\epsilon$.
It is direct to check that $V_\epsilon \in \M$ and its equilibrium
potentials are
$$
\mu_{V_\epsilon}=\lambda_\epsilon,\;\;
\nu_{V_\epsilon}=D_{b_\epsilon}\nu_{U_\epsilon}.
$$
Thus $V_\epsilon \in bi\MP$ and we can apply Theorem
\ref{potencialfU} to get that the matrix $V_\epsilon^{(\alpha)}$
posses a right equilibrium potential $\eta_\epsilon\ge 0$, that
is, for all $i$
$$
\sum\limits_{j} (V^{(\alpha)}_\epsilon)_{ij} (\eta_\epsilon)_j=1,
$$
which is equivalent to
$$
\sum\limits_{j}
\frac{(U_\epsilon)_{ij}^{\alpha}}{(b_\epsilon)_i^\alpha}
(\eta_\epsilon)_j=1.
$$
Hence
$$
U_\epsilon^{(\alpha)} \eta_\epsilon=b_\epsilon^{(\alpha)}.
$$
Recall that the matrix $U^{(\alpha)}$ is nonsingular. Since
obviously $U^{(\alpha)}_\epsilon \to U^{(\alpha)}$ as $\epsilon
\to 0$, we get
$$
\eta_\epsilon \to \eta=(U^{(\alpha)})^{-1} b^{(\alpha)},
$$
and the result follows. $\Box$
\end{proof}

\medskip
\begin{proof}({\bf Theorem \ref{power}.}) Consider first the case
where $U \in bi\MP$. We already know that $U^{(\alpha)}$ is
nonsingular and that it has left and right nonnegative equilibrium
potentials. Therefore, in order to prove that $U^{(\alpha)} \in
bi\MP$ is enough to prove that $(U^{(\alpha)})^{-1}$ is a $Z$-matrix
that is for $i\neq j$ we have
$((U^{(\alpha)})^{-1})_{ij}$ is non positive. An argument based on
permutations shows that it is enough to prove the claim for $i=1,
j=n$, where $n$ is the size of $U$.

Decompose $U^{(\alpha)}$ and its inverse as follows
$$
U^{(\alpha)}=\begin{pmatrix} A^{(\alpha)} &b^{(\alpha)} \cr
(c^{(\alpha)})' &d^\alpha
\end{pmatrix} \hbox{ and }
(U^{(\alpha)})^{-1}=\begin{pmatrix} \Omega &-\beta \cr -\alpha' &\delta
\end{pmatrix}
$$
We need to show that $\beta \ge 0$.
We notice that $\delta=\frac{det(A^{(\alpha)})}{det(U^{(\alpha)})}>0$
and that
$$
-A^{(\alpha)} \beta + \delta b^{(\alpha)}=0,
$$
which implies that
$$
b^{(\alpha)}=A^{(\alpha)}\left(\frac{\beta}{\delta}\right).
$$
Therefore, $\frac{\beta}{\delta}=\eta\ge 0$, where $\eta$ is the vector
given in Lemma \ref{generalequilibrium}. Thus $\beta \ge 0$ and the
result is proven  for the case $U\in bi\MP$.

Now consider $U$ the inverse of the $M$-matrix $M$. Using Theorem
2.5.3 in \cite{HJ}, we get the existence of two positive diagonal
matrices $D, E$ such that $DME$ is a strictly row and column
diagonally dominant $M$-matrix. Thus $V=E^{-1}UD^{-1}$ is in $bi\MP$
from where it follows that $V^{(\alpha)} \in bi\MP$. Hence,
$U^{(\alpha)}=E^{(\alpha)}V^{(\alpha)}D^{(\alpha)} $ is the
inverse of an $M$-matrix. The rest of the result is proven in a
similar way. $\Box$
\end{proof}

\medskip

\begin{proof}({\bf Theorem \ref{markovchain}.}) By hypothesis we have $U=\II-P$, where
$P\ge 0$ and $P\1\le \1$. We notice that $U$ is diagonally dominant on each row, that is
for all $i,j$
$$
U_{ii}\ge U_{ji}.
$$
Also we notice that $U=\II+PU$ and therefore $U_{ii}\ge 1$.

According to Theorem \ref{power}
we know that $H=(U^{(\alpha)})^{-1}$ is a row diagonally dominant $M$-matrix. The only thing
left to prove is that the diagonal elements of $H$ are dominated by one, that is
$H_{ii}\le 1$ for all $i$. It is enough to prove this for $i=n$, where $n$ is the size of $U$.

Consider the following decompositions
$$
\begin{array}{l}
U=\begin{pmatrix} A &b \cr c' &d
\end{pmatrix}\quad
U^{-1}=\begin{pmatrix}  \Lambda &-\omega \cr -\eta' &\gamma
\end{pmatrix}\quad
(U^{(\alpha)})^{-1}=\begin{pmatrix} \Omega &-\beta \cr -\alpha' &\delta
\end{pmatrix}
\cr
\cr
U^{-1}U^{(\alpha)}=\begin{pmatrix} \Xi &-\zeta \cr -\chi' &\rho
\end{pmatrix}
\end{array}
$$
A direct computation gives that
$$
\gamma=\rho \delta +\chi'\beta\ge \rho \delta.
$$
Since by hypothesis $\gamma \le 1$,
to conclude that $\delta \le 1$ it
is enough to prove that $\rho\ge 1$.
We have
$$
\rho=(1-p_{nn})U^{\alpha}_{nn}-\sum\limits_{j\neq n} p_{nj}U^{\alpha}_{jn}
=U^{\alpha}_{nn}-\sum\limits_{j} p_{nj}U^{\alpha}_{jn}=
U^{\alpha}_{nn}-\sum\limits_{j} p_{nj} U_{jn} U^{\alpha-1}_{jn}.
$$
On the other hand we have $U^{\alpha-1}_{jn}\le U^{\alpha-1}_{nn}$ and
$\sum\limits_{j} p_{nj} U_{jn}=U_{nn}-1$ from where we deduce
$$
\rho\ge U^{\alpha-1}_{nn}\ge 1.
$$
This finishes the case when $P\1\le \1$. The rest of the result is proven by using
$U'$ instead of $U$.
$\Box$
\end{proof}

\bigskip

\section{Proof of Theorem \ref{establecreciente->GUM}}
\medskip
Let $n$ be the dimension of $U$. Notice that $U$ is a GUM if and
only if $n\le 2$ or every principal sub-matrix of size $3$ is a
GUM.

Since by hypothesis the matrix $\II+tU$ is a bi-potential, it is
diagonally dominant
$$
1+tU_{ii}\ge tU_{ij},
$$
and we deduce that $U_{ii}\ge U_{ij}$ (by taking $t\to
\infty$). This proves the result when $n\le 2$. So, in the sequel
we assume $n\ge 3$.

Consider $A$ any principal sub-matrix of $U$, of size $3\times 3$.
Since $\II+t f(A)$ is a principal sub-matrix of $\II+t f(U)$, we
deduce that $\II+t f(A) \in bi\MP$ (as long as $\II+t f(U) \in
bi\MP$). If the result holds for the $3\times 3$ matrices, we
deduce that $A$ is GUM implying that $U$ is also a GUM.

\medskip

Thus, in what follows we consider that $U$ is a $3\times 3$ matrix
that verifies the hypothesis of the Theorem. After a suitable
permutation we can further assume that
$$
U=\begin{pmatrix} a&b_1 &b_2 \cr c_1 &d &\alpha \cr c_2 &\beta &e
\end{pmatrix},
$$
where $\alpha=\min\{U_{ij}: \; i\neq j\}=\min\{U\}$ and
$\beta=\min\{U_{ji}:\; U_{ij}=\alpha,\; i\neq j\}$.

\medskip

Since $U$ is diagonally dominant we have $\min\{a,d,e\}\ge
\alpha$. Take $f$ increasing such that $f(\alpha)=0$ and
$f(x)>0$ for $x>\alpha$. Then,
$$
\II+f(U)=\begin{pmatrix} 1+f(a)& f(b_1) & f(b_2) \cr f(c_1)
&1+f(d) &0 \cr f(c_2) &f(\beta) &1+f(e),
\end{pmatrix},
$$
is a $bi\MP$-matrix whose inverse we denote by
$$
\begin{pmatrix} \;\;\;\delta &-\rho_1 &-\rho_2 \cr -\theta_1
&\;\;\;\gamma_1 &-\gamma_2 \cr -\theta_2 &-\gamma_3
&\;\;\;\gamma_4
\end{pmatrix}.
$$
In particular we obtain
$$
\begin{pmatrix} 1+f(d) &0 \cr f(\beta) &1+f(e)
\end{pmatrix}^{-1}=\begin{pmatrix} \;\;\;\gamma_1
&-\gamma_2 \cr -\gamma_3 &\;\;\; \gamma_4
\end{pmatrix}-\frac{1}{\delta}\begin{pmatrix}
\theta_1\cr \theta_2
\end{pmatrix} \begin{pmatrix} \rho_1 \cr  \rho_2
\end{pmatrix}',
$$
and we deduce that
\begin{equation}
\label{caso1} 0=\gamma_2=\theta_1\rho_2.
\end{equation}

\medskip

\begin{itemize}

\item{}{\bf Case $\rho_2=0$.} We deduce that $f(b_2)=0$, and
then
\begin{eqnarray}
\label{caserho2=0} b_2=\alpha, \;\hbox{ and } c_2\ge \beta,
\end{eqnarray}
where the last conclusion follows from the definition of $\beta$.
Therefore we have that
\begin{equation}
\label{Ucasorho2=0} U=\begin{pmatrix} a&b_1 &\alpha \cr c_1 &d
&\alpha \cr c_2 &\beta &e
\end{pmatrix},
\end{equation}
and we should prove that $U$ is GUM.

Now consider another increasing function $g$ such that
$g(\beta)=0$ and $g(x)>0$ for $x>\beta$. Then,
$$
\II+g(U)=\begin{pmatrix} 1+g(a)& g(b_1) & 0 \cr g(c_1) &1+g(d) &0
\cr g(c_2) &0 &1+g(e)
\end{pmatrix}.
$$
Its inverse is of the form
$$
\begin{pmatrix} \;\;\;\tilde\delta &-\tilde\rho_1 &0 \cr -\tilde\theta_1
&\;\;\;\tilde\gamma_1 &0 \cr -\tilde\theta_2 &-\tilde\gamma_3
&\;\;\;\tilde\gamma_4
\end{pmatrix}.
$$
As before we deduce that
$0=\tilde\gamma_3=\tilde\theta_2\tilde\rho_1$.

\medskip

\begin{itemize}

\item{}{\bf Subcase $\tilde\theta_2=0$.} We have $g(c_2)=0$ which
implies $c_2=\beta$. Thus, in this situation we have that $U$ is
$$
U=\begin{pmatrix} a& b_1 & \alpha \cr c_1 & d &\alpha \cr \beta
&\beta &e
\end{pmatrix}.
$$
By permuting rows and columns $1,2$, if necessary, we can assume
that $b_1\le c_1$. Consider the situation where $c_1<\beta$, of
course implicitly we should have $\alpha<\beta$. Under a suitable
increasing transformation $h$ we have
$$
\II+h(U)=\begin{pmatrix} 1+h(a)& 0 & 0 \cr 0 &1+h(d) &0 \cr
h(\beta) &h(\beta) &1+h(e)
\end{pmatrix},
$$
and its inverse
$$
\begin{pmatrix} \frac{1}{1+h(a)} &0 &0 \cr 0
&\frac{1}{1+h(d)}  &0 \cr -\frac{h(\beta)}{(1+h(a))(1+h(e))}
&-\frac{h(\beta)}{(1+h(d))(1+h(e))} &\frac{1}{1+h(e)}
\end{pmatrix}.
$$
The sum of the third row is then
$$
\frac{1}{1+h(e)}\left(1-h(\beta)\left(\frac{1}{1+h(a)}+\frac{1}{1+h(d)}\right)\right),
$$
and this quantity can be made negative by choosing an appropriate function
$h$. The idea is to make $h(\beta)\to \infty$ and
$$
\frac{h(\beta)}{\max\{h(a),h(d)\}}\to 1.
$$
Therefore, $c_1\ge \beta$ and $U$ is a GUM.

\medskip

\item{}{\bf Subcase $\tilde\rho_1=0$.} We have $g(b_1)=0$ and then
$b_1\le \beta$. Take again an increasing function, denoted by
$\ell$, such that
$$
\II+\ell(U)=\begin{pmatrix} 1+\ell(a)& 0 & 0 \cr \ell(c_1)
&1+\ell(d) &0 \cr \ell(c_2) &0 &1+\ell(e)
\end{pmatrix},
$$
and its inverse
$$
\begin{pmatrix} \frac{1}{1+\ell(a)} &0 &0 \cr -\frac{\ell(c_1)}{(1+\ell(a))(1+\ell(d))}
&\frac{1}{1+\ell(d)}  &0 \cr
-\frac{\ell(c_2)}{(1+\ell(a))(1+\ell(e))} &0 &\frac{1}{1+\ell(e)}
\end{pmatrix}.
$$
The sum of the first column is
$$
\frac{1}{1+\ell(a)}\left(1-\frac{\ell(c_1)}{(1+\ell(d))}-\frac{\ell(c_2)}{(1+\ell(e))}\right),
$$
which can be made negative by repeating a similar argument as
before, if both $c_1>\beta$ and $c_2>\beta$.
\smallskip

So, if $c_1>\beta$ we conclude that $c_2 \le \beta$, but we know
that $c_2\ge \beta$  (see (\ref{caserho2=0})) and we deduce that
$c_2=\beta$. Thus, $\alpha\le b_1\le \beta <c_1$ and
$$
U=\begin{pmatrix} a& b_1 & \alpha \cr c_1 & d &\alpha \cr \beta
&\beta &e
\end{pmatrix},
$$
which is a GUM.

\medskip

Therefore we can continue under the hypothesis $c_1\le \beta \le
c_2$.

\begin{itemize}

\item{}{\bf Subsubcase $b_1<\beta$.} Again we must have $\alpha<\beta$.
Under this conditions we have that $c_2>\alpha$. Using an
increasing function $k$ we get
$$
\II+k(U)=\begin{pmatrix} 1+k(a)& 0 & 0 \cr k(c_1) & 1+k(d) &0 \cr
k(c_2) &k(\beta) &1+k(e)
\end{pmatrix},
$$
and its inverse is
$$
\begin{pmatrix} \frac{1}{1+k(a)}& 0 & 0 \cr & & \cr -\frac{k(c_1)}{(1+k(a))(1+k(d))} & \frac{1}{1+k(d)} &0 \cr
& & \cr
-\frac{k(c_2)(1+k(d))-k(\beta)k(c_1)}{(1+k(a))(1+k(d))(1+k(e))}
&-\frac{k(\beta)}{(1+k(d))(1+k(e))} &\frac{1}{1+k(e)}
\end{pmatrix}.
$$
The sum of the third row is
$$
\frac{1}{(1+k(e))}\left(1-\frac{k(c_2)}{1+k(a)}+\frac{k(\beta)k(c_1)}{(1+k(a))(1+k(d))}
-\frac{k(\beta)}{1+k(d)}\right).
$$
If $c_1<\beta$ we can assume that $k(c_1)=0$, and we can make this
sum to be negative by choosing large $k$. Thus we must have
$c_1=\beta$, in which case the sum under study is proportional to
\begin{equation}
\label{suma1}
1-\frac{k(c_2)}{1+k(a)}+\frac{k(\beta)^2}{(1+k(a))(1+k(d))}
-\frac{k(\beta)}{1+k(d)}
\end{equation}
If $c_2=\beta$ then
$$
U=\begin{pmatrix} a& b_1 & \alpha \cr \beta & d &\alpha \cr \beta
&\beta &e
\end{pmatrix},
$$
is a GUM. So we must analyze the case where $c_2>\beta$ in
(\ref{suma1}). We will arrive to a contradiction by taking an
asymptotic as before. Consider a fixed number $\lambda \in (0,1)$.
Choose a family $(k_r)_{r\in \NN}$ such that as $r \to \infty$
$$
k_r(\beta)\to \infty,\; \frac{k_r(\beta)}{k_r(c_2)}\to \lambda,\;
\frac{k_r(c_2)}{k_r(a)}\to 1,\; \frac{k_r(d)}{k_r(a)}\to \phi,
$$
where $\phi=1$ if $d>\beta$, and $\phi=\lambda$ if $d=\beta$. The
asymptotic of (\ref{suma1}) is then
$$
1-1+\frac{\lambda^2}{\phi}-\frac{\lambda}{\phi}.
$$
This quantity is strictly negative for the two possible values of
$\phi$, which is a contradiction, and therefore $c_2=\beta$.

\bigskip

To finish with the {\bf Subcase $\tilde\rho_1=0$}, which will in
turn finish with {\bf Case $\rho_2=0$}, we consider

\item{}{\bf Subsubcase $b_1=\beta$.} We recall that we are
under the restrictions $c_1\le \beta \le c_2$ and
$$
U=\begin{pmatrix} a& \beta & \alpha \cr c_1 & d &\alpha \cr c_2
&\beta &e
\end{pmatrix}.
$$
Notice that if $c_2=\beta$ then $U$ is GUM. So for the rest of
this subcase we assume $c_2>\beta$. Also if $c_1=\alpha$ we can
permute $1$ and $2$ to get
$$
\Pi U\Pi'=\begin{pmatrix} d& \alpha & \alpha \cr \beta & a &\alpha
\cr \beta &c_2 &e
\end{pmatrix},
$$
which is also in NBF, and $U$ is a GUM. Thus we can assume that
$c_1>\alpha$, and again of course we have $\alpha<\beta$.

\smallskip

Take an increasing function $m$ such that
$$
\II+m(U)=\begin{pmatrix} 1+m(a)& m(\beta) & 0 \cr m(c_1) & 1+m(d)
&0 \cr m(c_2) &m(\beta) &1+m(e)
\end{pmatrix},
$$
We take the asymptotic under the following restrictions:
$$
\frac{m(\beta)}{m(a)}\to \lambda \in (0,1),\;
\frac{m(c_1)}{m(a)}\to \lambda,\; \frac{m(e)}{m(a)} \to 1,\;
\frac{m(c_2)}{m(a)}\to 1,\;  \frac{m(d)}{m(a)}\to \phi,
$$
where $\phi=1$ if $d>\beta$, and it is $\lambda$ if $d=\beta$. The
limiting matrix for $\frac{1}{m(a)}(\II+m(U))$ is
$$
V=\begin{pmatrix} 1& \lambda & 0 \cr \lambda & \phi &0 \cr 1
&\lambda &1
\end{pmatrix},
$$
whose determinant is $\Delta=\phi-\lambda^2>0$. Therefore $V$ must
be in $bi\MP$. On the other hand the inverse of $V$ is given by
$$
V^{-1}=\frac{1}{\Delta}\begin{pmatrix} \phi& -\lambda & 0 \cr
-\lambda & 1 &0 \cr -(\phi-\lambda^2) &0 &\phi-\lambda^2
\end{pmatrix},
$$
and the sum of the first column is
$$
\frac{\lambda^2-\lambda}{\Delta}<0,
$$
which is a contradiction.
\end{itemize}

\end{itemize}

This finishes with the subcase $\rho_2=0$ and we return to
(\ref{caso1}) to consider now the following case

\item{}{\bf Case $\theta_1=0$.} Under this condition we get
$c_1=\alpha$ and
$$
U=\begin{pmatrix} a& b_1 & b_2 \cr \alpha & d &\alpha \cr c_2
&\beta &e
\end{pmatrix}.
$$
Consider the transpose of $U$ and permute on it  $2$ and $3$, to
obtain the matrix
$$
\tilde U=\begin{pmatrix} a& c_2 & \alpha \cr b_2 & e &\alpha \cr
b_1 &\beta &d
\end{pmatrix},
$$
where now $b_1\ge \beta$. Clearly the matrix $\tilde U$ verifies
the hypothesis of the Theorem and has the shape of
(\ref{Ucasorho2=0}), that is we are in the "case $\rho_2=0$" which
we already know implies that $\tilde U$ is GUM. Therefore $U$
itself is GUM.
\end{itemize}
$\Box$

\bigskip

\section{Filtered Matrices, sufficient conditions for classes $bi\MP$ and $\mathcal{T}$}
\label{filtrada}
A class of matrices of our interest is the
class of filtered matrices, which turn to be a generalization of
GUM matrices. They were introduced as operators in \cite{dell1998} to generalize
the class of self adjoint operators whose spectral decomposition is written in terms
of conditional expectations (see for instance \cite{bouleau1989},
\cite{dart1988} and \cite{martinez1994}).

The basic tool
to construct these matrices are partitions of ${\cal
J}_n=\{1,\cdots,n\}$. The components of a partition $\p$ are
called atoms. We denote by $\mathop{\thicksim}\limits^\p $ the
equivalence relation induced by $\p$. Then $i,j$ are in the same
atom of $\p$ if and only if $i\mathop{\thicksim}\limits^\p j $.

\medskip

A partition $\p$ is coarser or equal than $\q$ if the
atoms of $\q$ are contained in the atoms of $\p$. We denote this
(partial) order relation by $\p \preccurlyeq \q$. For example in
${\cal J}_4$ we have $\p=\{ \{1,2\},\{3,4\}\} \<
\q=\{\{1\},\{2\},\{3,4\}\}$. The coarsest partition is the trivial
one ${\cal N}=\{{\cal J}_n\}$ and the finer one is the discrete
partition ${\cal F}=\{\{1\},\{2\},\cdots,\{n\}\}$.

\begin{definition}
A filtration is an strictly increasing sequence of comparable
partitions $\FF=\{\p_0 \prec \p_1 \prec \cdots \prec \p_k\}$. A
filtration in wide sense is an increasing sequence of comparable
partitions $\GG=\{\p_0\< \p_1\< \cdots \< \p_k\}$.
\end{definition}

The difference between these two concepts is that in the latter
repetition of partitions is allowed.

\begin{definition} A filtration $\FF =\{\p_0\prec \cdots \prec \p_k\}$
is called dyadic if each non-trivial atom of $\p_s$ is divided
into two atoms of $\p_{s+1}$.
\end{definition}

The following example is the simplest dyadic filtration
$$
\FF=\{{\cal N}\prec \{\{1\},\{2,\cdots,n\}\}\prec \cdots \prec
\{\{1\},\{2\},\cdots,\{i\},\{i+1,\cdots,n\}\}\prec \cdots \prec
{\cal F}\}.
$$

Each partition $\p$ induces an incidence matrix $F=:F(\p)$ given
by
$$
F_{ij}=\begin{cases}1 & \hbox{ if } i\mathop{\thicksim}\limits^\p
j \cr 0 &\hbox{ otherwise }
\end{cases}
$$
A vector $v\in \RR^n$ is said to be $\p$-measurable if $v$ is
constant on the atoms of $\p$, that is
$$
i\mathop{\thicksim}\limits^\p j \Rightarrow v_i=v_j\;.
$$
This can be expressed in terms of standard matrix operations as
$$
F(\p) v=D_{w_\p} v,
$$
where $w_\p=F(\p)\1$ is the vector of sizes of the atoms. Recall
that $D_z$ is the diagonal matrix associated to the  vector $z$.
The set of $\p$-measurable vectors is a linear subspace of
$\RR^n$.

\begin{definition}
A matrix $U$ is said to be {\bf filtered} if there exists a
filtration in wide sense $\GG=\{\q_0 \< \q_1 \< \cdots \< \q_l\}$,
vectors $\mathfrak{a}_1,\cdots,\mathfrak{a}_l,\;
\mathfrak{b}_1,\cdots,\mathfrak{b}_l$ with the restriction that
$\mathfrak{a}_s,\mathfrak{b}_s$ are $\q_{s+1}$-measurable, such
that
\begin{equation}
\label{eq1} U=\sum\limits_{s=0}^\ell D_{\mathfrak{a}_s} F(\q_s)
D_{\mathfrak{b}_s},
\end{equation}
\end{definition}
There is no loss of generality if we assume that $\q_0={\cal N}$
and $\q_\ell={\cal F}$, that is $F(\q_0)=\1 \1'$ and
$F(\q_\ell)=\II$, the identity matrix. Let us see that (\ref{eq1})
can be simply written in terms of a filtration. Indeed, notice
that if $\mathfrak{a}_s$ and $\mathfrak{b}_s$ are
$\q_s$-measurable then
$$
D_{\mathfrak{a}_s} F(\q_s)
D_{\mathfrak{b}_s}=D_{\mathfrak{a}_s}D_{\mathfrak{b}_s}F(\q_s)=D_{\mathfrak{a}_s\odot
\mathfrak{b}_s} F(\q_s),
$$
where the vector $\mathfrak{a}_s\odot \mathfrak{b}_s$ is the
Hadamard product of $\mathfrak{a}_s$ and $\mathfrak{b}_s$ which is
also $\q_s$-measurable. Hence a sum of terms of the form
$$
D_{\mathfrak{a_s}}F(\q_s)D_{\mathfrak{b_s}}+D_{\mathfrak{a}_{s+1}}F(\q_{s+1})D_{\mathfrak{b}_{s+1}}+\cdots+
D_{\mathfrak{a}_{s+r}}F(\q_{s+r})D_{\mathfrak{b}_{s+r}},
$$
with $\p=\q_s=\cdots=\q_{s+r}$, can be reduced to the sum of two
terms as
$$
D_{C}F(\p)+D_{\mathfrak{a}_{s+r}}F(\p)D_{\mathfrak{b}_{s+r}},
$$
where $C=\sum\limits_{k=0}^{r-1} \mathfrak{a}_{s+k}\odot
\mathfrak{b}_{s+k}$, which is $\p$-measurable. In this way representation
(\ref{eq1}) can be written as
\begin{equation}
\label{reduced} U=\sum\limits_{s=0}^k
D_{C_s}F(\p_s)+D_{\mathfrak{m}_s}F(\p_s)D_{\mathfrak{n}_s}
\end{equation}
where $\FF=\{\p_0 \prec \p_1 \prec \cdots \prec \p_k\}$ is a
filtration, ${\cal N}=\p_0, \; {\cal F}=\p_k$, $C_s$ is
$\p_s$-measurable, $\mathfrak{m}_s, \mathfrak{n}_s$ are
$\p_{s+1}$-measurable and $\mathfrak{m}_k=0$. We shall always
consider this reduced representation of (\ref{eq1}), and we shall
say that $U$ is {\bf filtered} with respect to the filtration
$\FF$.

\medskip

If all $\mathfrak{m}_s,\mathfrak{n}_s$ are $\p_s$-measurable then
(\ref{eq1}) reduces to the form
\begin{equation}
\label{symmetric} U=\sum\limits_{s=0}^k D_{C_s+\mathfrak{m}_s\odot
\mathfrak{n}_s} F(\p_s),
\end{equation}
and $U$ is a symmetric matrix.

\medskip

We are mainly interested in a decomposition like (\ref{reduced})
with the vectors $\mathfrak{m}_s, \mathfrak{n}_s$ having the
following special structure:
\begin{equation}
\label{especial} \mathfrak{m}_s=\Gamma_s \odot p_s,\;
\mathfrak{n}_s=q_s,
\end{equation}
where $\Gamma_s$ is $\p_s$-measurable and $\{p_s,q_s\}$ is a
$\p_{s+1}$-measurable partition, that is, they are
$\p_{s+1}$-measurable $\{0,1\}$-valued vectors with disjoint
support: $p_s\odot q_s=0$ and $p_s+q_s=1$. If this is the case we
say that $U$ is a Special Filtered Matrix ({\bf SFM})
\begin{equation}
\label{SFM} U=\sum\limits_{s=0}^k D_{C_s}\;F(\p_s) + D_{\Gamma_s}
D_{p_s}\; F(\p_s) \; D_{q_s}.
\end{equation}
Notice that $\Gamma_k=0$.

\bigskip

%
% relacion entre CBF y filtrada
%

It is not difficult to see that every CBF matrix is filtered. This
is done by induction. Assume that
$$
U=\begin{pmatrix} A &\alpha \1_p\1_{n-p}'\cr \beta \1_{n-p}\1_p' &
B
\end{pmatrix}.
$$
Define $\p_0={\cal N}$ and
$\p_1=\{\{1,\cdots,p\},\{p+1,\cdots,n\}\}$. Take $C_0=\alpha \1_n,
\Gamma_0=(\beta-\alpha)\1_n, p_0=(\0_p, \1_{n-p})',
q_0=(\1_p,\0_{n-p})'$. Then
$$
D_{C_0} F(\p_0)+D_{\Gamma_0}D_{p_0} F(\p_0) D_{q_0}=
\begin{pmatrix} \alpha \1_p\1_p' &\alpha \1_p\1_{n-p}'\cr \beta
\1_{n-p}\1_p' & \alpha \1_{n-p}\1_{n-p}'
\end{pmatrix}.
$$
The key step is that $A-\alpha,B-\alpha$ are also in CBF. We have
that $C_0, \Gamma_0$ are $\p_0$-measurable and $p_0,q_0$ is a
$\p_1$-measurable partition. We also notice that if $0\le
\alpha\le \beta$ then $C_0\ge 0, \Gamma_0\ge 0$.

\medskip

The induction also shows that $U$ can be decomposed as in
(\ref{SFM}) $U=\sum\limits_{s=0}^k D_{C_s}F(\p_s) + D_{\Gamma_s}
D_{p_s}F(\p_s)D_{q_s}$, where $\FF=\{\p_0 \prec \cdots \prec
\p_k\}$ is a dyadic filtration, $C_s, \Gamma_s$ are
$\p_s$-measurable and $\{p_s,q_s\}$ is a $\p_{s+1}$-measurable
partition.

\medskip

We summarize now the representation form for the class of CBF, NBF
and GUM matrices.

\medskip

\begin{proposition}
\label{CBFisSFM} $V$ is a permutation of a CBF if and only if
there exists a dyadic filtration $\FF=\{\p_0 \prec \cdots \prec
\p_k\}$, a sequence of vectors $C_0,\cdots,C_k, \;
\Gamma_0,\cdots, \Gamma_k$ verifying that for all $i$: $C_s,
\Gamma_s$ are $\p_{s}$-measurable and a sequence $\{p_s, q_s\}$ of
 $\p_{s+1}$-measurable partitions, such that
$$
V=\sum\limits_{s=0}^k D_{C_s}F(\p_s) + D_{\Gamma_s} D_{p_s}
F(\p_s)D_{q_s},
$$
that is $V$ is a SFM. Also $V$ is a permutation of an increasing
CBF matrix if and only if there is a decomposition where
$\Gamma_0, C_s,\Gamma_s:\; s=1\dots,k$, are nonnegative. On the top
of this $V$ is a nonnegative matrix if and only if $C_0$ is
nonnegative.

Moreover, $V$ is a GUM if and only if $\; C_s,\Gamma_s:\;
i=0,\cdots,k$  are non-negative, and for $i=0,\dots,k-1$ it holds
\begin{equation}
\label{americana} \Gamma_s\le C_{s+1}+\Gamma_{s+1}.
\end{equation}
Finally, $V$ is an ultrametric matrix if and only if there is a decomposition
with $\Gamma_s=0$ for all $s$. $\Box$
\end{proposition}

\begin{remark}
\label{diadico} We can assume without loss of generality that each
$p_s, q_s$ is obtained as follows. The nontrivial atoms ${\cal
A}_1,\cdots {\cal A}_r$ of $\p_s$ are divided into the new atoms
$$
{\cal A}_{1,1},{\cal A}_{1,2}\cdots {\cal A}_{r,1},{\cal
A}_{r,2}
$$
of $\p_{s+1}$. Consider also ${\cal B}_1,\cdots {\cal B}_r$ the
set of trivial atoms in $\p_s$ (that is the singleton atoms). Take
$q_s$ be the indicator of ${\cal A}_{1,1}\cup \cdots \cup {\cal
A}_{r,1}$, $p_s$ be the indicator of ${\cal A}_{1,2}\cup \cdots
\cup {\cal A}_{r,2}\cup {\cal B}_1\cup \cdots \cup {\cal B}_r$ and
$\Gamma_s=0$ on ${\cal B}={\cal B}_1\cup \cdots \cup {\cal B}_r$,
which is $\p_s$-measurable. We point out that the partition
$\p_{s+1}$ is obtained from $\p_s$ refined by $p_s$. Also the
important relation holds
\begin{equation}
\label{eqdiadica} D_{p_s}F(\p_s)p_s=D_{p_s}F(\p_{s+1})\1
\end{equation}

\end{remark}

\medskip

\begin{example} Consider the CBF matrix
$$
U=\begin{pmatrix} a &\alpha_2  &\alpha_1 &\alpha_1\cr \beta_2 &b
&\alpha_1 &\alpha_1 \cr \beta_1 &\beta_1 &c &\hat\alpha_2 \cr
\beta_1 &\beta_1 &\hat\beta_2 &d
\end{pmatrix}.
$$
$U$ is a NBF if the following constraints are verified:
$\alpha_1\le \beta_1,\; \alpha_1\le
min\{\alpha_2,\hat\alpha_2\},\; \beta_1\le
min\{\beta_2,\hat\beta_2\},\; \alpha_2\le \beta_2,\;
\hat\alpha_2\le \hat\beta_2$ and finally the diagonal dominates
over each row and column, that is $\beta_2\le min\{a,b\},\;
\hat\beta_2\le min\{c,d\}$.

$U$ is filtered with respect to the dyadic filtration
$\p_0=\{1,2,3,4\}\prec \p_1=\{\{1,2\},\{3,4\}\} \prec
\p_2=\{\{1\},\{2\},\{3\},\{4\}\}$ and can be written as
$$
U=D_{C_0} F(\p_0)+D_{\Gamma_0} D_{p_0} F(\p_0) D_{q_0}+D_{C_1}
F(\p_1)+D_{\Gamma_1} D_{p_1} F(\p_1) D_{q_1}+D_{C_2} F(\p_2),
$$
where
$$
C_0=\begin{pmatrix} \alpha_1\cr \alpha_1\cr \alpha_1\cr
\alpha_1\end{pmatrix},\; \Gamma_0=\begin{pmatrix}
\beta_1-\alpha_1\cr \beta_1-\alpha_1\cr \beta_1-\alpha_1\cr
\beta_1-\alpha_1\end{pmatrix},\; p_0=\begin{pmatrix} 0\cr 0\cr
1\cr 1\end{pmatrix},\; q_0=\begin{pmatrix} 1\cr 1\cr 0\cr
0\end{pmatrix}
$$
$$
C_1=\begin{pmatrix} \alpha_2-\alpha_1\cr \alpha_2-\alpha_1\cr
\hat\alpha_2-\alpha_1\cr \hat\alpha_2-\alpha_1\end{pmatrix},\;
\Gamma_1=\begin{pmatrix} \beta_2-\alpha_2\cr \beta_2-\alpha_2\cr
\hat\beta_2-\hat\alpha_2\cr
\hat\beta_2-\hat\alpha_2\end{pmatrix},\; p_1=\begin{pmatrix} 0\cr
1\cr 0\cr 1\end{pmatrix},\; q_1=\begin{pmatrix} 1\cr 0\cr 1\cr
0\end{pmatrix},
$$
and
$$
C_2=\begin{pmatrix} a-\alpha_2\cr b-\alpha_2\cr c-\hat\alpha_2\cr
d-\hat\alpha_2 \end{pmatrix}.
$$
The constrains are translated into: the positivity of these
vectors and the ones induced by (\ref{americana}). We point out
that we can also choose, for example,
$\Gamma_1=(0,\beta_2-\alpha_2,0,\hat\beta_2-\hat\alpha_2)'$, but
in this case $\Gamma_1$ is not $\p_1$-measurable.  As we will see
in subsection (\ref{sec-algorithm}) this measurability condition
will play an important role.
\end{example}

\begin{example} Consider the nonnegative CBF matrix

$$
U=\begin{pmatrix} 2 &2  &2 \cr 2 &2 &1 \cr 2 &1 &2
\end{pmatrix}.
$$
This matrix is a SFM and can be decomposed as in (\ref{SFM}).
Nevertheless no such decomposition can have all terms nonnegative.
In particular no permutation of $U$ is an increasing CBF.
\end{example}

\bigskip

\begin{remark} Notice that the class of CBF is stable under Hadamard functions.
Nevertheless there are examples of filtered matrices for which
$f(U)$ is not filtered. Consider the matrix
$$
U=D_\alpha F_1+ D_{\mathfrak{a}} F_1 D_{\mathfrak{b}} + D_\beta
F_2,
$$
where $F_1=F({\cal N})=\1\1'$ and $F_2=\II$. We have $\alpha$ is a
constant vector, and we confound it with the constant $\alpha\in
\RR $. The vectors $\mathfrak{a},\; \mathfrak{b}, \beta$ are all
${\cal F}$-measurable. Then $U$ is filtered and moreover
\begin{equation}
\label{eq3} U=\alpha+\mathfrak{a}\mathfrak{b}'+D_\beta
\end{equation}

Take $\alpha=\beta=0$ and $\mathfrak{a}=(2,3,5,7)'$ and
$\mathfrak{b}=(11,13,17,19)'$. Then all the entries of $U$ are
different. As $f$ runs all possible functions $f(U)$ runs over all
$4\times 4$ matrices. This implies that some of them can not be
written as in (\ref{eq3}), because in this representation we have
at most 13 free variables. Still is possible that each $f(U)$ is
decomposable as in (\ref{eq1}) using maybe a different filtration.
A more detailed analysis shows that this is not the case. For
example if we choose the filtration ${\cal N}\prec
\{\{1,2\},\{3,4\}\} \prec {\cal F}$ then every matrix $V$ filtered
with respect to this filtration verifies that
$V_{13}=V_{23}=V_{14}=V_{24}$.
\end{remark}

\bigskip

Matrices of the type $F(\p)$ are related to conditional
expectations (in probability theory). Indeed, let
$\p=\{{\cal A}_1,{\cal A}_2,\cdots,{\cal A}_r\}$ and
$n_\ell=\#({\cal A}_\ell)$ be the size of each atom. It is direct
that $w=w_\p=F(\p)\1$ is a $\p$-measurable vector that verifies
$w_i=n_\ell$  for $i \in {\cal A}_\ell$. Then
$$
\EE_\p=D_w^{-1}\;
F(\p)=F(\p)\; D_w^{-1}
$$
is the matrix of conditional expectation with respect to the
$\sigma$-algebra generated by $\p$. This matrix $\EE=\EE_\p$
satisfies:
$$
\begin{array}{ll}
&\EE \EE=\EE; \; \EE'=\EE; \; \EE\1=\1;\cr &\forall v\; \EE v
\hbox{ is } \p-\hbox{measurable};\cr  &\hbox{ if } v \hbox{ is }
\p-\hbox{measurable, then } \EE v=v.
\end{array}
$$
Therefore, $\EE$ is the orthogonal projection over the subspace of
all $\p$-measurable vectors. In the case of the trivial partition
${\cal N}$ one gets $\EE_{\cal N}=\frac{1}{n} \1\1'$ the mean
operator.

\bigskip

\begin{remark} \label{rm1} The $L^2$ space associated to
$\{1,\cdots,n\}$ endowed with the counting measure, is identified with
$\RR^n$ with the standard euclidian scalar product.
In this way each vector of $\RR^n$ can be seen as a
function in $L^2$, and $\EE$ is an orthogonal projection.
The product $D_v \EE$ (as matrices) is the product of the operators
$D_v$ and $\EE$, where $D_v$ is the multiplication by
the function $v$. Notice that $\EE D_v$ and $\EE(v)$ are
quite different. The former is an operator (a matrix) and the
latter is a function (vector). They are related by $\EE(v) = \EE
D_v (\1)$, where $\1$ is the constant function.
\end{remark}

\medskip

Let $\p,\;\q$ be two partitions, then
$\p\<\q$ is equivalent to
$\EE_\p\EE_\q=\EE_\q\EE_\p=\EE_\p$. This commutation relation can
be written as a commutation relation for $F(\p)$ and $F(\q)$.
In fact,
$$
\begin{array}{ll}
F(\p)F(\q)&=\EE_\p  D_{w_\p} \EE_\q D_{w_\q}= \EE_\p \EE_\q
D_{w_\p}  D_{w_\q}=\cr &= \EE_\p  D_{w_\p}  D_{w_\q}=F(\p)
D_{w_\q};\cr F(\q)F(\p)&=(F(\p)F(\q))'=D_{w_\q} F(\p).
\end{array}
$$

\medskip

\subsection{An algorithm for filtered matrices: conditions to be in $bi\MP$}
\label{sec-algorithm}

In this section we explain a backward algorithm to determine when
a filtered matrix is in class $bi\MP$. Assume that $U$ has a
representation as in (\ref{eq1})
$$
U=\sum\limits_{s=0}^\ell D_{\mathfrak{a}_s} F(\q_s)
D_{\mathfrak{b}_s},
$$
where we assume further that $\mathfrak{a}_s, \mathfrak{b}_s$ are
all non-negative. In particular $U$ is a non-negative matrix.

\medskip

We introduce the conditional expectations
$\EE_s=\EE_{\q_s}=D_{F(\q_s)\1}^{-1}\; F(\q_s)$
and the normalized factors: $a_s=\mathfrak{a}_s\odot F(\q_s)\1,\;
b_s=\mathfrak{b}_s$. Then $U$ can be written as
$$
U=\sum\limits_{s=0}^\ell D_{a_s} \EE_s
D_{b_s}=\sum\limits_{s=0}^\ell a_s \EE_s b_s,
$$
where we have identified vectors (functions) and the operator of
multiplication they induce. We shall use this notation throughout
this section. Finally, we remind that $\EE_\ell=\II$.

\medskip

We can now use the algorithm developed in \cite{dell1998} to study
the inverse of $\II+U$. In what follows, we take the convention
$0\cdot\infty=0/0=0$. This algorithm is defined by the backward
recursion starting with the values
$\lambda_\ell=\mu_\ell=\kappa_\ell=1,\; \sigma_\ell=(1+a_\ell
b_\ell)^{-1}$ and for $s=\ell-1,\cdots,0$:
\begin{equation}
\label{algoritmo}
\begin{matrix}
&\lambda_{s}&=&\lambda_{s+1}[1-\sigma_{s+1}a_{s+1}\EE_{s+1}(\kappa_{s+1}b_{s+1})];\cr
&\mu_{s}&=&\mu_{s+1}[1-\sigma_{s+1}b_{s+1}\EE_{s+1}(\kappa_{s+1}a_{s+1})];\cr
&\kappa_s&=&\EE_{s+1}(\lambda_s)=\EE_{s+1}(\mu_s); \cr
&\sigma_s&=&(1+\EE_s(\kappa_sa_sb_s))^{-1},
\end{matrix}
\end{equation}
from where we obtain the recursion
\begin{equation}
\label{recursion}
\kappa_{s-1}=\EE_s(\kappa_s)-\frac{\EE_s(\kappa_sa_s)
\EE_s(\kappa_s b_s)}{1+\EE_s(\kappa_s a_s b_s)}.
\end{equation}
The algorithm continues until some $\lambda$ or $\mu$ are negative
otherwise we arrive to $s=0$. If this is the case then $\II+U$ is
nonsingular and its inverse is of the form $\II-N$ where
$$
N = \sum_{s=0}^\ell \sigma_s\lambda_sa_s\EE_sb_s\mu_s.
$$
We also have that
$$
\lambda_{-1}=(\II-N)\1, \hbox{ and } \mu_{-1}=(\II-N)'\1,
$$
where $\lambda_{-1}, \mu_{-1}$ are obtained from the first two
formulae in (\ref{algoritmo}) for $s=-1$. Therefore, if they are
also non-negative the matrix $\II+U$ is a $bi\MP$-matrix.

\medskip

In this way we have that a sufficient condition for $\II+U$ to be
a $bi\MP$-matrix, is that the algorithm works for $s=\ell,\cdots,0$
and all the $\lambda, \mu$ are nonnegative, including
$\lambda_{-1}, \mu_{-1}$. In this situation we have that $\lambda$
(and $\mu$) is a decreasing non-negative sequence of vectors.

Sufficient treatable conditions  involve the recurrence
(\ref{recursion}). Starting from $\kappa_\ell=1$ we assume this
recurrence  has a solution such that $\kappa_s \in [0,1]$ for all
$s=\ell,\cdots,-1$. We shall study closely this recursion for the
class of SFM, and we shall obtain sufficient conditions to have
$\II+U$ in $bi\MP$.

\medskip

Before studying this problem, we discuss further the algorithm. We
have the following relations:
$$
\begin{array}{l}
\left(\II+\sum\limits_{k=s}^\ell a_k\EE_k
b_k\right)^{-1}\!\!\!=\II-\sum\limits_{k=s}^\ell \sigma_k
\lambda_k a_k \EE_k b_k \mu_k=\II-N_s, \cr \cr
\lambda_{s-1}=(\II-N_s)\1,\; \mu_{s-1}=(\II-N_s)'\1
\end{array}
$$
That is, the algorithm imposes that all the matrices:
$$
\II+a_\ell \EE_\ell b_\ell,\cdots, \II+\sum\limits_{k=s}^\ell
a_k\EE_k b_k, \cdots, \II+\sum\limits_{k=0}^\ell a_k\EE_k
b_k=\II+U,
$$
are in class $bi\MP$.

\bigskip

We now assume that $U$ is a SFM with a decomposition like
$$
U=\sum\limits_{s=0}^k D_{C_s}\;F(\p_s) + D_{\Gamma_s} D_{p_s}\;
F(\p_s) \; D_{q_s},
$$
where $\FF=\p_0\prec \cdots \prec \p_k$ is a filtration,
$C_s, \Gamma_s$ are nonnegative $\p_s$-measurable and
$\{p_s, q_s\}$ is a $\p_{s+1}$-measurable partition. Again
we put $\EE_s=D_{F(\p_s)\1}^{-1}\; F(\p_s)$ and the normalized factors:
$$
c_s=C_s\odot F(\p_s)\1,\; \gamma_s=\Gamma_s \odot F(\p_s)\1,
$$
which are $\p_s$-measurable. Since diagonal matrices commute we
get that $U$ has a representation of the form
$$
U=\sum\limits_{s=0}^k c_s\EE_s+ \gamma_s p_s \EE_s q_s,
$$
with $\gamma_k=0$. In the previous algorithm we can make two steps
at each time and consider  $\kappa_s$ in place of $\kappa_{2s}$,
$\lambda_s$ instead of $\lambda_{2s+1}$, $l_s$ instead of
$\lambda_{2s}$. We also introduce $d_s=1/\kappa_s$ to simplify
certain formulae (this vector can take the value $\infty$). We
get, starting from $\kappa_k=l_k=1, \sigma_k=(1+c_k)^{-1}$, that
for $s=k-1,\cdots,0$
$$
\begin{array}{lrll}
&\lambda_{s} & = & \sigma_{s+1}l_{s+1};\cr
&l_s&=&\lambda_{s}[1-\gamma_{s}p_{s}
\EE_{s}(q_{s}/(c_{s+1}+d_{s+1}))]; \cr &\kappa_s &=&\EE_s(l_s);\cr
&\sigma_s & =&1/(1+\kappa_sc_s)=d_s/(c_s+d_s).
\end{array}
$$
Similar recursions hold for $\mu, m$, which are the analogous of
$\lambda, l$. Relation (\ref{recursion}) takes the form
\begin{equation}
\label{iteracion}
\frac{1}{d_s}=\EE_s\left(\frac{1}{c_{s+1}+d_{s+1}}\right)-
\gamma_s\,\EE_s\left(\frac{p_s}{c_{s+1}+d_{s+1}}\right)
\EE_s\left(\frac{q_s}{c_{s+1}+d_{s+1}}\right)
\end{equation}

The inverse of $\II+U$ is $\II-N$ where
\begin{equation}
\label{N} N=\sum_{s=0}^k c_s \sigma_s l_s \EE_s m_s+
\sum_{s=0}^{k-1} \gamma_s \lambda_s p_s \EE_s q_s \mu_s
=\sum_{s=0}^k c_s \sigma_s l_s \EE_s m_s+ \gamma_s \lambda_s p_s
\EE_s q_s \mu_s.
\end{equation}
Again $\lambda_{-1}=(\II-N)\1=s_0l_0$, and similarly
$\mu_{-1}=s_0m_0$.

\medskip

Let us introduce the following function
$$
\rho_s=\EE_s(p_s)p_s+\EE_s(q_s)q_s.
$$

\begin{theorem}
\label{invM} Assume that the backward recursion (\ref{iteracion})
has a non-negative solution starting with $d_k=1$. Assume moreover
that this solution verifies for $s=k-1,\cdots,0$
\begin{equation}
\label{desigualdad} \rho_s \gamma_s\le c_{s+1}+d_{s+1}.
\end{equation}
Then $\lambda_s,l_s,\mu_s,m_s,s_s :\; i=k,\cdots,0$, as well as
$\lambda_{-1}, \mu_{-1}$, are well defined and nonnegative.
Therefore, $\II+U \in bi\MP$ and its inverse is $\II-N$ where $N$ is
given by (\ref{N}).
\end{theorem}

\medskip

The proof of this result is based on the following lemma.

\begin{lemma}
\label{fundamental} Assume $x, y $ are nonnegative vectors and
$\;\EE$ is a conditional expectation. If $\;x \EE(y)\le 1$ then
$\EE(xy)\le 1$.
\end{lemma}

\begin{proof}
We first assume that $y$ is strictly positive. Since $x\le
1/\EE(y)$ and $\EE$ is an increasing operator, we have
$$
\EE(x y)\le \EE\left(\frac{1}{\EE(y)}
y\right)=\frac{\EE(y)}{\EE(y)}=1.
$$
For the general case consider $(y+\epsilon\1)/(1+\epsilon
|x|_\infty)$ instead of $y$ and pass to the limit $\epsilon\to 0$.
$\Box$
\end{proof}

\smallskip

\begin{proof}{({\bf Theorem \ref{invM}})}
We notice that condition (\ref{desigualdad}) implies that
$$
\frac{ q_s \gamma_s}{c_{s+1}+d_{s+1}}\EE_s(q_s)\le 1.
$$
Since $\gamma_s$ is $\EE_s$-measurable and $q_s=q_s^2$ we obtain
$$
\gamma_s
\EE_s\left(\frac{q_s}{c_{s+1}+d_{s+1}}\right)=\EE_s\left(\frac{\gamma_s
q_s^2}{c_{s+1}+d_{s+1}}\right).
$$
This last quantity is bounded by one by Lemma \ref{fundamental}.
Similarly we have
$$
\gamma_s \EE_s\left(\frac{p_s}{c_{s+1}+d_{s+1}}\right)\le 1,
$$
which implies that the algorithm is not stopped, all the
coefficients are non-negative including $\lambda_{-1}, \mu_{-1}$.
$\Box$
\end{proof}

\bigskip

\begin{corollary}
\label{gum} Assume that for $s=k-1,\cdots,0$ we have
\begin{equation}
\label{eq.gum} \rho_s \gamma_s\le c_{s+1}+\gamma_{s+1}.
\end{equation}
Then the recursion (\ref{iteracion}) has a nonnegative solution
that verifies (\ref{desigualdad}). In particular, $\II+tU$ is in
class $bi\MP$ for all $t\ge 0$, and $U$ is in $bi\MP$ if it is
nonsingular.
\end{corollary}

\begin{proof}
Let us consider first the case $t=1$. We prove by induction that
$\gamma_s\le d_s$. For $i=k$ we have $0=\gamma_k\le d_k=1$. We
point out that if we multiply in (\ref{iteracion}) by $\gamma_s$
we get
$$
\frac{\gamma_s}{d_s}=\EE_s\left(\frac{\gamma_s}{c_{s+1}+d_{s+1}}\right)-
\EE_s\left(\frac{\gamma_s p_s }{c_{s+1}+d_{s+1}}\right)
\EE_s\left(\frac{\gamma_s q_s}{c_{s+1}+d_{s+1}}\right),
$$
which is of the form $x+y-xy$, where $x=\EE_s\left(\frac{\gamma_s
p_s }{c_{s+1}+d_{s+1}}\right)$. The inequality (\ref{eq.gum}), the
induction hypothesis $\gamma_{s+1}\le d_{s+1}$ and Lemma
\ref{fundamental} imply $0\le x\le 1, 0\le y\le 1$. In particular
$$
0\le \frac{\gamma_s}{d_s}\le 1,
$$
and the induction is completed. Theorem \ref{invM} shows that
$\II+U$ is in class $bi\MP$. We notice that $tU$ also verifies
condition (\ref{eq.gum}) because this condition is homogeneous,
and the result follows. $\Box$
\end{proof}

\begin{remark} We notice that condition (\ref{eq.gum}) can be
expressed in terms of the original coefficients $C, \Gamma$ in the
dyadic case. In fact (see (\ref{eqdiadica}))
$$
p_s\EE_s(p_s)=D_{p_s} D_{F(\p_s)\1}^{-1}\;
F(\p_s)p_s=D_{p_s}D^{-1}_{F(\p_s)\1} F(\p_{s+1})\1,
$$
which implies that
$$
\rho_s=(1/F(\p_s)\1)\odot F(\p_{s+1})\1.
$$
Then inequality (\ref{eq.gum}) is
$$
\Gamma_s\le C_{s+1}+\Gamma_{s+1},
$$
which is the condition for having a GUM (see (\ref{americana})) .
We mention here that condition (\ref{eq.gum}) is more general than
having a GUM, as the following example shows.
\end{remark}

\begin{remark}
Consider the matrix $U_\beta$
$$
U_\beta=\begin{pmatrix} 1 &0 &0 &0\cr 0 &1 &0 &0\cr \beta &\beta
&1 &0\cr \beta &\beta &0 &1 \end{pmatrix}=D_{\Gamma_0} D_{p_0}
F(\p_0) D_{q_0}+ \II,
$$
where $\p_0={\cal N}$, $\Gamma_0=\beta(1,1,1,1)'\le
C_1=(1,1,1,1)'$. We compute $c_0=0,\; \gamma_0=4\beta, \;
c_1=C_1,\; \gamma_1=0$ and also $\rho_0=1/2$.

It is direct to check that $U_\beta^{-1}=U_{-\beta}$. Then for all
$\beta\ge 0$ the matrix $U_\beta \in \M$. Also $U_\beta \in bi\MP$
if and only if $0\le \beta \le 1/2$. When $\beta\ge 0$ the
condition (\ref{americana}), that is
$$
\Gamma_0\le C_1+\Gamma_1,
$$
is equivalent to $ \beta \le 1$. Then, this condition does not
ensure that $U\in bi\MP$ (this happens because the filtration is not dyadic).
Nevertheless, the analogue condition in terms
of the normalized factors (\ref{eq.gum})
$$
\rho_0 \gamma_0 \le c_1 + \gamma_1,
$$
which is equivalent to $\beta \le 1/2$, is the right one.
\end{remark}

\begin{corollary} Assume that
\begin{equation}
\label{frances} \rho_s\gamma_s\leq \sum_{r=s+1}^k c_r,
\end{equation}
hold for $s=k-1,\cdots,0$. Then the recursion (\ref{iteracion})
has a nonnegative solution that verifies (\ref{desigualdad}). In
particular, $\II+tU$ is in class $bi\MP$ for all $t\ge 0$, and $U$
is in $bi\MP$ if it is nonsingular.
\end{corollary}

\begin{proof} Consider the set of inequalities
$$
\rho_s\gamma_s\vee\xi_s\leq c_{s+1}+\xi_{s+1},
$$
for $i=k-1,\cdots,0$.  A non-negative solution is given by
$\xi_s=\sup\{0,\; \gamma_0\rho_0-\sum_{r=1}^s
c_r,\cdots,\gamma_k\rho_k-\sum_{r=k+1}^s
c_r,\cdots,\gamma_{s-1}\rho_{s-1}-c_s\}$.  The hypothesis is that $\xi_k=0$.
Moreover we have that $\xi_s$ is
$\p_s$-measurable.

We show, using a backward recursion that $\xi_s\le d_s$. Indeed,
$1/\xi_s=\EE_s(1/\xi_s) \ge (c_{s+1}+\xi_{s+1})^{-1} $ by
construction while $1/d_s\leq \EE_s((c_{s+1}+d_{s+1})^{-1})$. Then
the inequality $\rho_s\gamma_s\leq c_{s+1}+\xi_{s+1} $ implies
$\rho_s\gamma_s\leq c_{s+1}+d_{s+1}$, from where the result holds
(see Theorem \ref{invM}). $\Box$
\end{proof}

\bigskip

\subsection { Conditions for class $\mathcal{T}$ and proof of Theorem \ref{CBFisTau} }

\bigskip

\begin{theorem}
\label{filteredistau}
Assume that $U$ has a decomposition
$$
U=\sum\limits_{s=0}^\ell a_s\;  \EE_s \; b_s,
$$
where $a_s,\; b_s$ are nonnegative $\EE_{s+1}$-measurable. Then
$U$ belongs to the class $\mathcal{T}$ and moreover
$$
\tau(U)=\inf\{t>0:\; (\II+tU)^{-1}\1 \ngtr 0 \hbox{ or }
\1'(\II+tU)^{-1} \ngtr 0\}
$$
In particular if $\tau(U)<\infty$ then  $\II+\tau(U)\, U\in bi\MP$.

\end{theorem}

\begin{remark}
Since the set of nonsingular matrices is open, then in the
previous result when $\tau(U)<\infty$, we have  for $t>\tau(U)$
sufficiently close to $\tau(U)$, that the matrix $\II+t\,U$ is
nonsingular.
\end{remark}

\medskip

Theorem \ref{filteredistau} states that every filtered matrix,
with a nonnegative decomposition, is in class $\mathcal{T}$ which
proves Theorem \ref{CBFisTau}.

\bigskip

\begin{proof}({\bf Theorem \ref{filteredistau}.})
A warning about the use of vectors and functions. Here we consider
vectors or functions on $\{1,\dots,n\}$ indistinctively. Thus for
two vectors $a,b$ the product $a b$ makes sense as the product of
two functions, which corresponds to the Hadamard product of the
vectors. Also an expression as $(1+ab)^{-1}$ is the vector whose
components are the reciprocals of the components of $1+ab$. We
also recall that $(a)_i$ is the $i$-th  component of $a$.

First, for $p=0,\dots, \ell$ consider the matrices
$$
U(p)=\sum\limits_{s=p}^\ell a_s\;  \EE_s \; b_s.
$$
We notice that $U(0)=U$. We shall prove that $\tau_p=\tau(U(p))$
is increasing in $p$ and $\tau_\ell=\infty$.

We rewrite the algorithm for $\II+t U$. This takes the form
$\lambda_\ell(t)=\mu_\ell(t)=\kappa_\ell(t)=1,\;
\sigma_\ell(t)=(1+t\,a_\ell b_\ell)^{-1}$ and for
$p=\ell-1,\cdots,0$:
\begin{equation}
\label{algoritmo-t}
\begin{matrix}
&\lambda_{p}(t)&=&\lambda_{p+1}(t)[1-\sigma_{p+1}(t)\,t\,
a_{p+1}\EE_{p+1}(\kappa_{p+1}(t)b_{p+1})];\cr
&\mu_{p}(t)&=&\mu_{p+1}(t)[1-\sigma_{p+1}(t)\,t\,
b_{p+1}\EE_{p+1}(\kappa_{p+1}(t)a_{p+1})];\cr
&\kappa_p(t)&=&\EE_{p+1}(\lambda_p(t))=\EE_{p+1}(\mu_p(t)); \cr
&\sigma_p(t)&=&(1+\EE_p(\kappa_p(t)ta_pb_p))^{-1},
\end{matrix}
\end{equation}
Also $\lambda_{-1}(t), \; \mu_{-1}(t)$ are defined similarly. If
$\lambda_s(t),\; \mu_s(t),\; \sigma_s(t):\; s=\ell,\dots, p$ are
well defined then
$$
(\II+t U(p))^{-1}=\II-N(p,t),
$$
where
\begin{equation}
\label{Npt} N(p,t)=\sum\limits_{s=p}^\ell \sigma_s(t)
\lambda_s(t)\, t \, a_s\EE_s b_s \mu_s(t).
\end{equation}
If $\lambda_s(t),\; \mu_s(t),\; \sigma_s(t):\; s=\ell,\dots, p$
are nonnegative then $N(p,t)\ge 0$, and $(\II+t U(p))\in \M$.
Moreover, $\lambda_{p-1}(t)$ and $\mu_{p-1}(t)$ are the right and
left equilibrium potentials of $(\II+t U(p))$
$$
(\II+t U(p))\lambda_{p-1}(t)=\1, \hbox{ and } \mu'_{p-1}(t)(\II+t
U(p))=\1'.
$$
So, if they are nonnegative, we have $\II+tU(p) \in bi\MP$. In
particular we have that
$$
(\II+ta_\ell\,\EE_\ell\, b_\ell)^{-1}=(\II+tU(\ell))^{-1}=\II-t
(1+t\,a_\ell b_\ell)^{-1} a_\ell \,\EE_\ell \,b_\ell.
$$
Since $\EE_\ell=\II$ we obtain that
$\lambda_{\ell-1}=\mu_{\ell-1}=(1+t\,a_\ell b_\ell)^{-1}$. This
means that $\II+tU(\ell) \in bi\MP$ for all $t\ge 0$. Therefore
$\tau_\ell=\infty$ and the result is true for $U(\ell)$. This
implies in particular that $\tau_{\ell-1}\le \tau_\ell$. Assume
the following inductive hypothesis
\begin{itemize}
\item{} $\tau_{p+1}\le \dots \le \tau_\ell$ ;

and for $q=p+1,\dots,\ell$
\item{} $\tau_q\!=\!\inf\{t>0\!: \lambda_{q-1}(t) \ngeqq 0 \hbox{ or }
\mu_{q-1}(t) \ngeqq 0\}\!=\!\inf\{t>0\!: \lambda_{q-1}(t) \ngtr 0
\hbox{ or } \mu_{q-1}(t) \ngtr 0\}$;
\item{} $\lambda_s(t), \mu_s(t)$, for
$s=\ell,\dots,q-1$,
are strictly positive for  $t\in [0,\tau_q)$ ;
\item{} If $\tau_q<\infty$ we have  $\II+\tau_q
U(q) \in bi\MP$.
\end{itemize}

\medskip

The case $\tau_{p+1}=\infty$ is clear. Indeed, fix $t\ge 0$. From
Lemma \ref{M->Mt}, $\II+t U(p+1) \in bi\MP$ and its equilibrium
potential are strictly positive, that is $\lambda_p(t)>0,
\mu_p(t)>0$. Thus, $\II+t U(p)$ is nonsingular, its inverse is
$\II-N(p,t)$, where $N(p,t)\ge 0$ is given by (\ref{Npt}). Hence,
$\II+tU(p)\in \M$. We conclude that
$$
\tau_p=\inf\{t>0: \II+t U(p) \notin bi\MP\}=\inf\{t>0:
\lambda_{p-1}(t)\ngeqq 0 \hbox{ or }  \mu_{p-1}(t)\ngeqq 0 \}.
$$
So, if $\tau_p=\infty$ we have, from Lemma \ref{M->Mt} that
$$
\lambda_{p-1}(t)>0,\; \mu_{p-1}(t)> 0,
$$
and the induction step holds in this case.

Now if $\tau_p<\infty$, by continuity we have $\II+\tau_p U(p) \in
bi\MP$, and we shall prove later on that $\lambda_{p-1}(t),
\mu_{p-1}(t)$ are strictly positive in $[0,\tau_p)$.

\bigskip

We analyze now the case $\tau_{p+1}<\infty$. We first notice that
in the algorithm the only possible problem is with the definition
of $\sigma_p(t)$. Since $\sigma_p(\tau_{p+1})>0$ the algorithm is
well defined, by continuity, for steps $\ell,\dots,p$  on an
interval $[0,\tau_{p+1}+\epsilon]$, for small enough $\epsilon$.
This proves that the matrix $\II+t U(p)$ is nonsingular in that
interval, and that $\lambda_{p-1}, \mu_{p-1}$ exist in the same
interval.

\bigskip

Now, for a sequence $t_n\downarrow \tau_{p+1}$ either
$\lambda_p(t_n)$ or $\mu_p(t_n)$ has a negative component. Since
there are a finite number of components we can assume with no lost
of generality that for a fixed component $i$ we have
$(\lambda_p(t_n))_i<0$. Then, by continuity we get that
$(\lambda_p(\tau_{p+1}))_i=0$, which implies (by the algorithm)
that $(\lambda_{p-1}(\tau_{p+1}))_i=0$.

%If there exists another sequence $s_n\downarrow \tau_{p+1}$, such
%that $(\lambda_p(s_n))_i\ge 0$ we get a contradiction. In fact the
%function $z\to (\lambda_p(z))_i$ is a rational function of $z$ and
%therefore is analytic in a (complex) neighborhood of $\tau_{p+1}$.
%Since it is not the constant $0$, such sequence $(s_n)$ cannot
%exist. This implies that for a small $\delta>0$ all the components
%of $\lambda_p(t)$ that are not positive they must be negative on
%$t\in (\tau_{p+1},\tau_{p+1}+\delta)$.

\medskip

Assume now that for some $t>\tau_{p+1}$ the matrix $\II+tU(p) \in
bi\MP$. Again by Lemma \ref{M->Mt} we will have that $\II+\tau_{p+1}
U(p) \in bi\MP$ but its equilibrium potential will satisfy
$\lambda_{p-1}(\tau_{p+1})>0$, which is a contradiction. Therefore
we conclude that $\tau_p\le \tau_{p+1}$.

\bigskip

The conclusion of this discussion is that the matrix $\II+tU(p)$,
for $t\in [0,\tau_{p+1}]$, is nonsingular and its inverse is
$\II-N(p,t)$, with $N(p,t)\ge 0$. That is $\II+tU(p)\in \M$ and
therefore
$$
\tau_p=\inf\{t>0: \II+t U(p) \notin bi\MP\}=\inf\{t>0:
\lambda_{p-1}(t)\ngeqq 0 \hbox{ or }  \mu_{p-1}(t)\ngeqq 0 \},
$$
and by continuity $\II+\tau_p U(p) \in bi\MP$.

To finish the proof we need to show that $\tau_p$ coincides with
$$
S=\inf\{t>0: \lambda_{p-1}(t)\ngtr 0 \hbox{ or } \mu_{p-1}(t)\ngtr
0 \}.
$$
It is clear that $S\le \tau_p$. If $S< \tau_p$ then,  due to
Lemma \ref{M->Mt}, we have that both $\lambda_{p-1}(S)>0$ and
$\mu_{p-1}(S)> 0$, which is a contradiction and then $S=\tau_p$.
This shows that $\lambda_{p-1}(t), \mu_{p-1}(t)$ are strictly
positive for $t\in [0,\tau_p)$, and the induction is proven.
$\Box$
\end{proof}

\bigskip

\begin{remark}
It is possible to prove that $\kappa_p(\tau_p)>0$ when
$\tau_p<\infty$, but this is not central to our discussion.
\end{remark}


\bigskip

\begin{thebibliography}{99}

\bibitem{ando1980} T. Ando, Inequalities for $M$-matrices.
{\it Linear and Multilinear Algebra,}
{\bf 8} (1980), 291-316.

\bibitem{neummann2005} R.B. Bapat, M. Catral, M. Neummann,
On functions that preserve $M$-matrices and inverse $M$-matrices.
{\it Linear and Multilinear Algebra ,} {\bf 53} (2005), 193-201.


\bibitem{bouleau1989} N.~Bouleau, Autour de la variance comme
forme de Dirichlet, {\it S\'eminaire de Th\'eorie du Potentiel
{8}, Lecture Notes in Mathematics} {1235} (1989), 39--53.

\bibitem{dart1988} P.~Dartnell, S.~Mart\'{\i}nez and
J.~San Mart\'{\i}n, Op\'erateurs filtr\'es et cha\^{\i}nes de
tribus invariantes sur un espace probabilis\'e d\'enombrable. {\it
S\'eminaire de Probabilit\'es {XXII} Lecture Notes in Mathematics}
{1321}, (1988), Springer-Verlag.

\bibitem{dell1996} C.~Dellacherie, S.~Mart\'{\i}nez and
J.~San Mart\'{\i}n,  Ultrametric matrices and induced Markov
chains. {\it Advances in Applied Mathematics,} {\bf 17} (1996),
169--183.

\bibitem{dell2000} C.~Dellacherie, S.~Mart\'{\i}nez and
J.~San Mart\'{\i}n, Description of the Sub-Markov kernel associated to
generalized ultrametric matrices. An algorithmic approach. {\it
Linear Algebra and its Applications,} {\bf 318} (2000), 1-21.


\bibitem{dell1998} C.~Dellacherie, S.~Mart\'{\i}nez,
J.~San Mart\'{\i}n and D. Ta\"{\i}bi. Noyaux potentiels associ\'es
\`a une filtration. {\it Ann. Inst. Henri Poincar\'e Prob. et
Stat.} {\bf 34} (1998), 707--725.

\bibitem{fiedler1983} M. Fiedler, H. Schneider, Analytic Functions of $M$-matrices and generalizations.
{\it Linear and Multilinear Algebra,} {\bf 13} (1983), 185-201.

\bibitem{HJ} R. Horn and C. Johnson. {\it Topics in matrix analysis}.
Cambridge University Press (1991).

\bibitem{martinez1994} S.~Mart\'{\i}nez,  G.~Michon and
J.~San Mart\'{\i}n, Inverses of ultrametric matrices are of
Stieltjes types, {\it SIAM J. Matrix Analysis and its
Applications,} {\bf 15} (1994), 98--106.

\bibitem{mcdonald1995} J.J.~McDonald, M.~Neumann, H.~Schneider and
M.J.~Tsatsomeros, Inverse $M-$matrix inequalities and generalized
ultrametric matrices. {\it Linear Algebra and its Applications,}
{\bf 220} (1995), 321--341.

\bibitem{micchelli1979} C.A. Micchelli, R.A. Willoughby, On functions which preserve Stieltjes
matrices. {\it Linear Algebra and its Applications,} {\bf 23} (1979), 141-156.

\bibitem{nabben1994} R.~Nabben and R.S.~Varga,  A linear algebra
proof that the inverse of a strictly ultrametric matrix is a
strictly diagonally dominant Stieljes matrix,  {\it SIAM J. Matrix
Analysis and its Applications,} {\bf 15} (1994), 107--113.

\bibitem{nabben1995} R.~Nabben and R.S.~Varga, Generalized
ultrametric matrices -- a class of inverse $M$-matrices. {\it
Linear Algebra and its Applications,} {\bf 220} (1995), 365--390.

\bibitem{neummann1998} M. Neummann, A conjectured concerning the
Hadamard product of inverses of $M$-matrices. {\it Linear Algebra
and its Applications,} {\bf 285} (1998), 277-290.

\bibitem{varga1968}  R.S. Varga, Nonnegatively posed problems and completely monotonic
functions. {\it Linar Algebra and its Applications,} {\bf 1} (1968), 329-347.


\end{thebibliography}
\end{document}